\documentclass[11pt]{article}
\usepackage{graphicx}
\usepackage{amsmath}
\usepackage{amssymb}
\usepackage{amscd}
\usepackage{amsthm}
\usepackage{hyperref}
\usepackage{color}
\usepackage{enumerate}
\usepackage{caption}
\usepackage{subcaption}
\usepackage{adjustbox}
\usepackage[top=1in, bottom=1in, left=1in, right=1in]{geometry}
\usepackage[square,sort,comma,numbers]{natbib}
\usepackage[]{algorithm2e}
\usepackage{placeins}
\usepackage{authblk}

\usepackage{changes}

\title{Detecting transient rate-tipping using Steklov averages and Lyapunov vectors}
\date{\today}
\author[1]{Alanna Hoyer-Leitzel}
\author[2]{Alice Nadeau}
\author[3]{Andrew Roberts}
\author[4]{Andrew Steyer}
\affil[1]{Mt. Holyoke College}
\affil[2]{University of Minnesota}
\affil[3]{Cerner Corporation}
\affil[4]{Sandia National Laboratory}

\begin{document}

\maketitle

\begin{abstract}
A wide variety of physical systems ranging from the firing of neurons to eutrophication of lakes to the presence of Arctic summer sea ice exhibit a phenomenon known as tipping. In mathematical models, tipping can be caused by bifurcations, noise, and the rate at which parameters are changing in time \cite{ashwin2012}. Because traditional methods in dynamical systems are usually concerned with the long-term behavior of the system, these methods are not always able to detect the transient dynamics characteristic of rate-tipping. In this paper, we consider one- and two-dimensional dynamical systems with nonautonomous parameters that exhibit rate-tipping, as defined as not tracking the evolution of stable equilibria (QSEs) in the corresponding autonomous systems. We find that nonautonomous stability spectra in the form of Steklov averages and their derivatives appear to be correlated with transient rate-tipping in systems with unique QSEs or with parameters that change at a constant rate. Furthermore, for systems in two dimensions and higher, comparison of the angle between leading Lyapunov vectors of different trajectories admits a possible criterion for detecting rate-tipping. Our heuristic results add to the body of work dedicated to studying and understanding the phenomenon of rate-tipping. 
\end{abstract}

\section{Introduction}

The discussion of tipping points and critical transitions in climate, ecology, and related fields is widespread and increasingly urgent. In this discussion, it has become clear that the traditional mathematical models of a critical transition from bifurcation and catastrophe theories cannot account for some of the behaviors found in these systems (see \cite{lenton2013} and the references therein). For example, in autonomous bifurcation theory and in some forms of nonautonomous bifurcations, a control parameter is considered static and fixed at all time, and the question asked for applications is about the sensitivity of the parameter: what magnitude of error could there be in the estimation of a parameter before the system exhibits a critically different behavior? However, this discussion on tipping points and critical transitions posits that control parameters of a system are not static, and evolve over time. Examples of systems that exhibit rate-tipping (also called rate-induced tipping), show that the rate at which a parameter is changing can cause critical transitions in the dynamics of a system.

A time-dependent parameter recasts a model from an autonomous differential equations model to a nonautonomous model. In this paper we use methods and techniques for the approximation of stability spectra for nonautonomous differential equations as possible indicators of tipping. We find a correspondence of these tools with tipping, and even the possibility that they may be able to provide earlier detection of tipping than some previously proposed indicators. 

This paper consists of the following: In the rest of this section, we give background and definitions for rate-tipping and stability spectra for nonautonomous systems. Section 2 contains the methods used to numerically investigate tipping, and Section 3 contains examples and results. Section 4 is a discussion.

\subsection{Rate-tipping: Definitions and Context}
 
\begin{figure}[t]
        \begin{subfigure}[b]{.23\textwidth}
        \centering
        \includegraphics[width=1\textwidth]{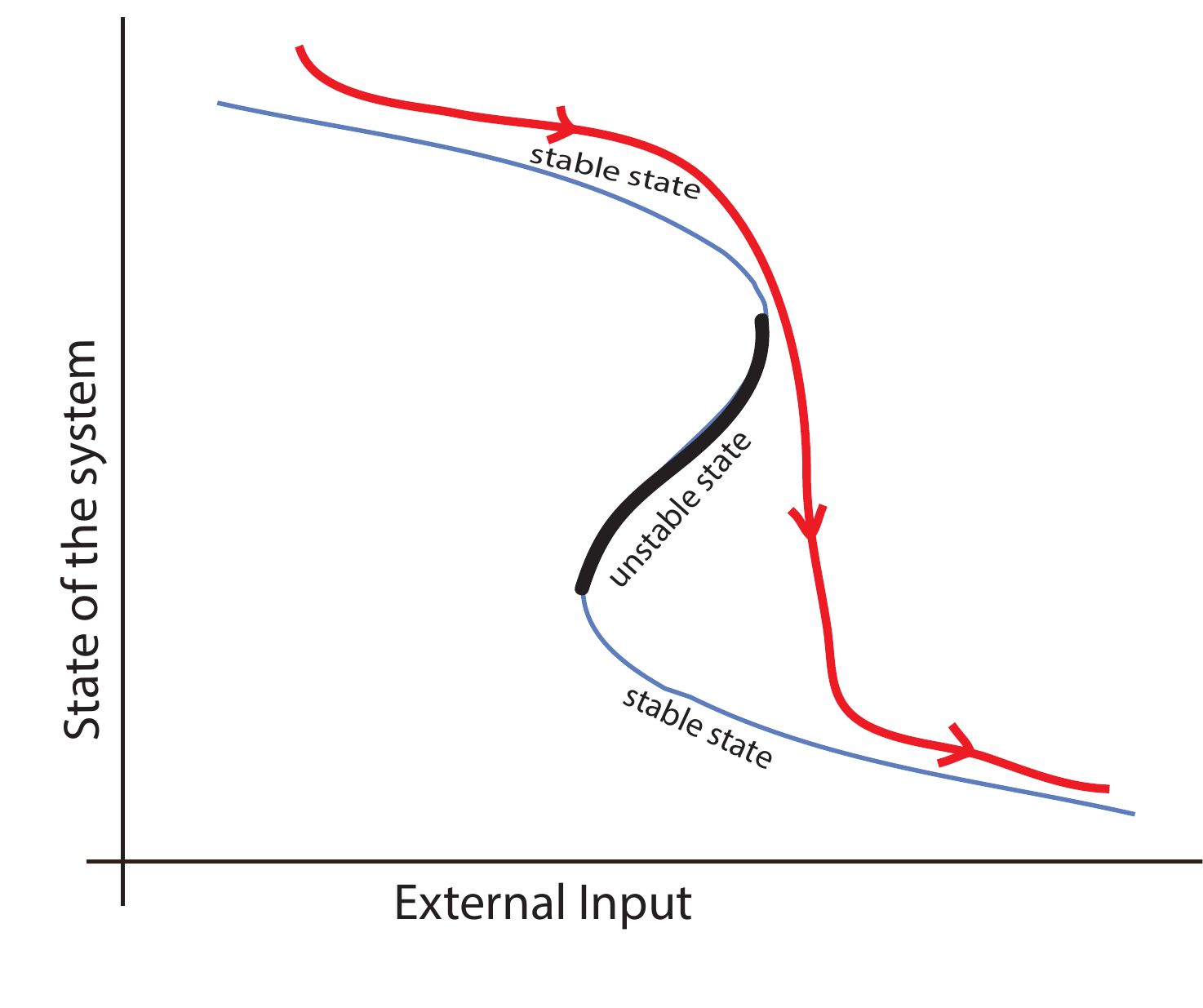}
        \caption{Bifurcation tipping}
        \label{Fig-BifTipping}
        \end{subfigure}
        \begin{subfigure}[b]{.23\textwidth}
        \vspace{22pt}
        \includegraphics[width=1\textwidth]{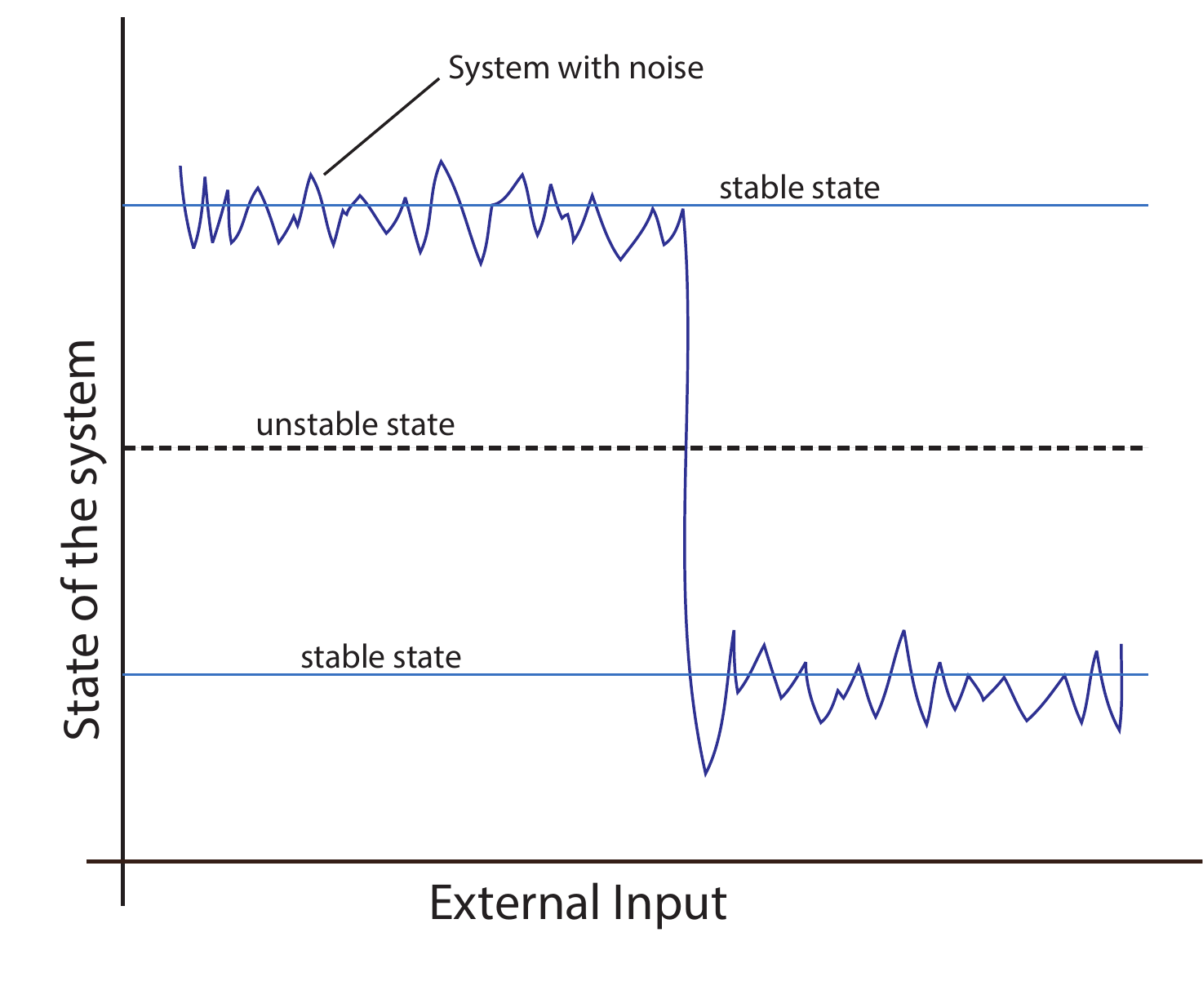}
        \caption{Noise tipping }
        \label{Fig-NoiseTipping}
        \end{subfigure}
        \begin{subfigure}[b]{.23\textwidth}
        \includegraphics[width=1\textwidth]{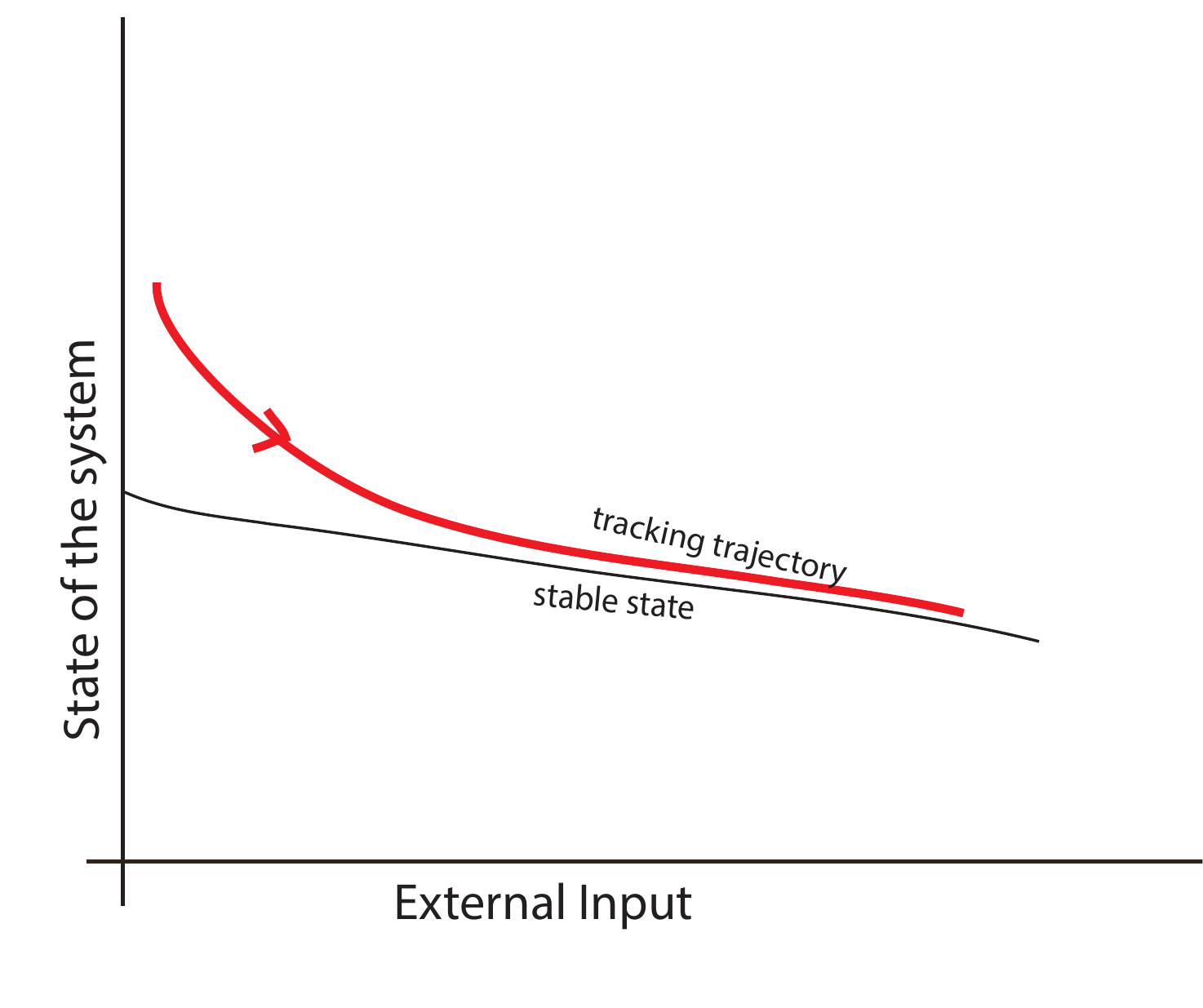} \caption{System tracking }
        \label{Fig-rate1}
        \end{subfigure}
        \begin{subfigure}[b]{.23\textwidth}
        \includegraphics[width=1\textwidth]{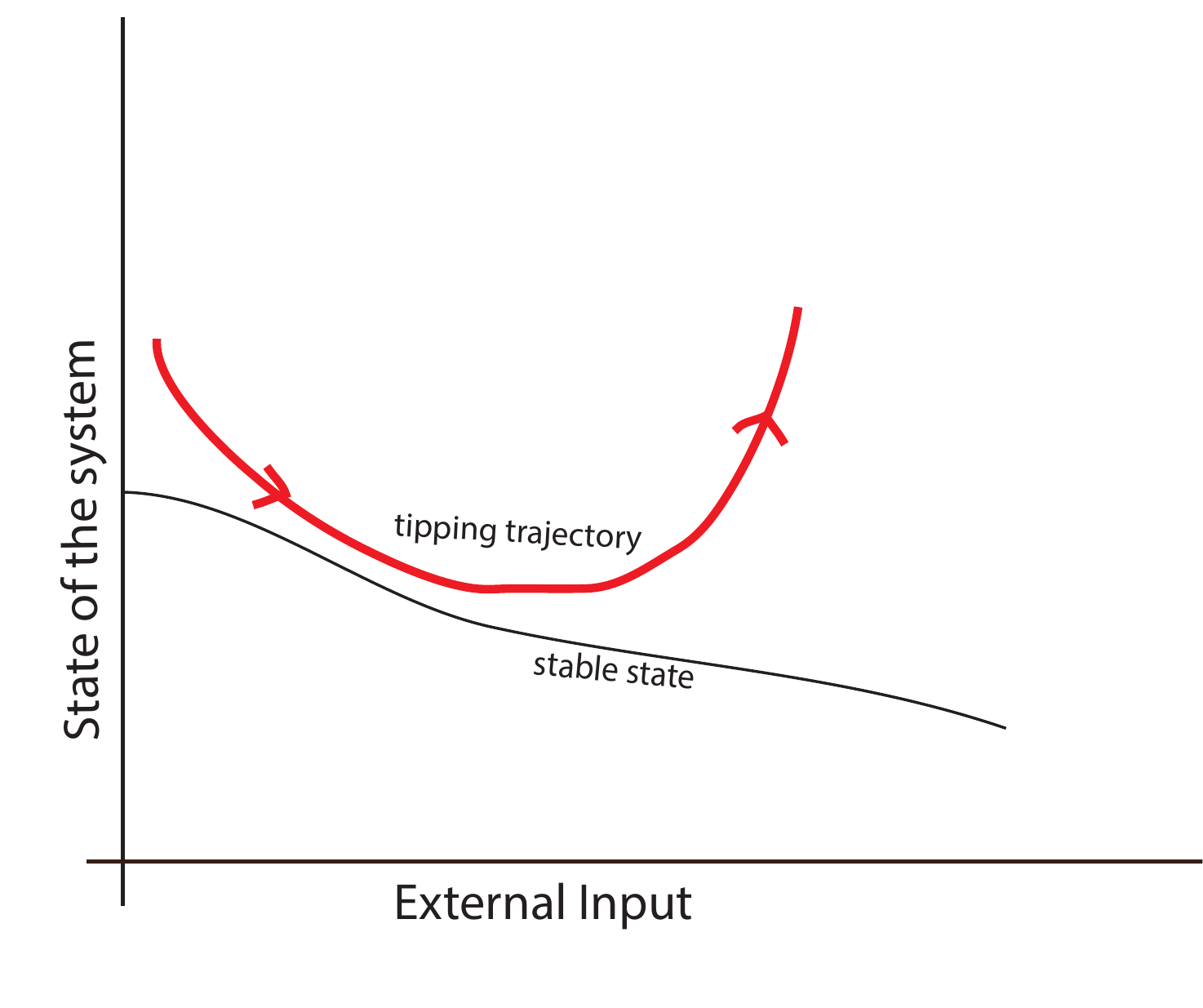} 
        \caption{Rate tipping }
        \label{Fig-rate2}
        \end{subfigure}
        \caption{Examples of Bifurcation, Noise, and Rate-tipping given by Ashwin et al. \cite{ashwin2012}.  In (a) the trajectory tips due to a bifurcation in a parameter, in (b) the trajectory tips due to noise in the system. Figures (c) and (d) illustrate rate-tipping. In (c) the trajectory tracks a stable quasi-static equilibrium until a critical rate in the change of the stable state, at which point it can no longer track the stable quasi-static equilibrium and tips away from the stable state (d).}
        \label{Fig-BifNoiseTipping-Examples}
    \end{figure}

 Tipping points are characterized by a sudden, qualitative shift in the behavior or state of the system due to a relatively small change in inputs \cite{lenton2008,scheffer2009} (e.g., a lake which becomes eutrophic due to increased nutrient run-off). Often a system that has tipped is difficult---or even impossible---to return back to its original state, making the study of predicting and preventing this tipping phenomenon highly important.  Ashwin et al (2012) \cite{ashwin2012} classified tipping phenomena into three categories based on the underlying mechanism which causes the system to tip. They find that tipping may result from 
\begin{itemize}
\item \textit{bifurcations} in state space (see Figure \ref{Fig-BifTipping}),
\item \textit{noise}, where noise within the system causes a change in state (see Figure \ref{Fig-NoiseTipping}), and 
\item \textit{rate}, where the rate at which a parameter is changing causes a change in state (see Figures \ref{Fig-rate1},\ref{Fig-rate2}).
\end{itemize}

In this paper we focus on rate-induced tipping. We look at examples of systems of the form $x \in \mathbb{R}^n$, 
\begin{equation}
\dot{x}=f(x,\lambda(t))
\label{eq:ode}
\end{equation}
where $\lambda(t)$ is a function of $t$ which is as smooth as the system, at least $C^1$. We define a \textit{quasi-static equilibrium}, or \textit{QSE}, as a set of equilibrium solutions (usually a smooth curve) of the associated autonomous systems where $\lambda(t) \equiv \lambda$ is constant. A QSE is called stable or unstable, depending on the stability of the equilibrium solution $\tilde{x}(\lambda)$ in the associated autonomous system. A trajectory of the nonautonomous system $\dot{x}=f(x,\lambda(t))$ is said to \textit{track} a stable QSE if $|x(t)-\tilde{x}(\lambda)|<R$ for some critical radius $R>0$. 

The dynamics of the nonautonomous system \eqref{eq:ode} can be drastically different than those of its associated autonomous systems. For example, in a delayed bifurcation, like that studied by Baer and Rinzel in \cite{baer1989}, the passage of $\lambda(t)$ through a traditional Hopf bifurcation value in the associated autonomous systems may allow the system to maintain the pre-bifurcation value dynamics, but only when the rate of $\lambda(t)$ is `slow enough'. In this paper and in most rate-tipping examples, we ignore possible effects of bifurcations by examining systems where the values of $\lambda(t)$ do not cross a traditional bifurcation value. In this case, the different dynamics occur when trajectories fail to track a stable QSE that is changing at a critical rate \cite{ashwin2012}. In fact, a stable QSE of the nonautonomous system may not attract solutions asymptotically (as it would in the set of autonomous systems). For example, in the nonautonomous saddle-node example in \cite{ashwin2012} at non-tipping rates of $\lambda(t)$, there is a forward attracting invariant line of the system that tracks the QSE. (See \cite{KS2005} for a definition of forward attracting invariant sets in nonautonomous systems.) 
\textit{In essence, the study of rate-tipping is about quantifying the boundary between changes in the time-dependent parameter that maintain the dynamics of the associated autonomous systems and changes in parameters that are too fast to maintain the autonomous dynamics, particularly when no bifurcation is present.}

Rate tipping has been investigated using various mathematical theories. Ashwin, Wieczorek, et al's first paper on rate tipping \cite{ashwin2012} characterizes a trajectory as \textit{tipping} when it fails to track a QSE, and analytically derives a relationship between a critical tipping rate and the tracking radius, assuming a linear rate of change of $\lambda(t)$, called \textit{steady drift} (which is also referred to as \textit{linear ramping} in \cite{wieczorek2011}).  In \cite{ashwin2012} they also consider rate-tipping in multi-scale fast-slow systems -- work that is carried forward by Perryman and Wieczork in \cite{perryman2014} using examples of parameter shifts that are asymptotically constant. In a third paper by Ashwin, Wieczorek, and Perryman \cite{ashwin2015}, rate-tipping problems are placed in the realm of nonautonomous systems theory, and tipping is characterized using pullback attractors (see \cite{KS2005} for a definition of pullback attractors) by showing that invariant pullback attractors can track stable QSE, and thus the trajectories attracted to them track the stable QSE as well. These pullback attractors are generalizations of the invariant lines in the saddle-node normal form example of \cite{ashwin2012}, whose existence implies there are trajectories that do not tip. 

The negative definition for tipping that a trajectory does not tip if it tracks a stable QSE lacks a robustness for other situations. In \cite{ashwin2015}, other characterizations of what it means for a trajectory to `not tip' are put forward. One such characterization is \textit{end-point tracking}, where where a trajectory may not track a QSE for some radius $R>0$ for all times $t$ but does asymptotically approach a QSE in the limit as $t\to\pm\infty$. This definition implies that the trajectory eventually tracks a stable QSE, and can be applied to systems with asymptotically constant parameters. However, both negative definitions of tipping are dependent on the distance of a trajectory from a stable QSE, but does not encompass all possible critical transitions that can be found in a system with a rate-dependent parameter. We find in this paper that tracking can be a reasonable characterization of `not-tipping' in cases where the time-dependent parameter is a coordinate shift, but not in other cases, like the multiplicative rate in the resource-consumer model example \cite{siteur2016ecosystems}. 

Recently, people have studied the interactions of rate- and noise-tipping in the same system. Ritchie and Sieber have a series of papers looking at this combination of tipping mechanisms, \cite{ritchie2016early,ritchie2016dynamic,ritchie2016probability}. In \cite{ritchie2016early}, the authors consider systems with asymptotically constant parameter drift and interpret rate-tipping as the breaking of a heteroclinic connection, which is also done analytically in \cite{perryman2015thesis}. Rate-tipping and noise have also been combined in applications as in \cite{chen2015} wherein the authors investigate systems that exhibit spatially periodic patterns through the lens of rate and noise tipping.  

One motivational application for rate-tipping is the Compost Bomb problem \cite{wieczorek2011}, where a steep enough linear increase in parameter leads an increase of 100$^{\circ}$C on Earth's temperature in a period of about 10 years. This problem had direct influence on the first paper on rate-tipping \cite{ashwin2012}. However, the Compost Bomb system is quite stiff and is very difficult to numerically integrate accurately. Instead, as an application of the techniques in this paper, we consider the ecological resource-consumer model given by \cite{siteur2016ecosystems}. This two dimensional system also exhibits rate tipping due to a linear parameter change but is less stiff than the Compost Bomb.  We are able to detect tipping in the Resource-Consumer model with our numerical methods, and these results are given in Section \ref{Section-Examples}. 

The lack of true equilibrium solutions in the nonautonomous system $\dot{x}=f(x,\lambda(t))$ and the fact that a QSE may not be a trajectory necessitates different techniques for analysis than stability theory for equilibrium solutions in autonomous systems. In this paper, we compare numerical methods for stability in nonautonomous systems as indicators of tipping and the critical thresholds developed in the previous study of rate-tipping phenomena. 

\subsection{Stability in Nonautonomous Systems: Context and Theory}
The stability of a solution to a differential equation is usually defined by the long term asymptotic behavior. Rate-tipping involves finding the critical rate at which there is a loss of stability locally in time, and hence rate-tipping is a question about the transient dynamics of a system. In this section, we introduce the spectra used to analyze stability of a trajectory in a nonautonomous system and their numerical approximations. We show that the numerical approximations of these spectra can be used to find the short term growth and decay rates of a trajectory, hence their usefulness for determining if a trajectory tips.

We characterize the stability of a time-dependent solution $x(t;x_0,t_0)$ of $\dot{x} = f(x,\lambda(t))$ through the stability properties of the zero solution of the associated linear variational equation 
\begin{equation}
\dot{u} = A(t)u, \quad A(t) = D_x f(x(t;x_0,t_0),\lambda(t))
\label{eq:lve}
\end{equation}
where $D_x$ denotes the partial derivative of $f(x,\lambda)$ with respect to $x$.  We remark that the stability properties of the zero solution $u(t) \equiv 0$ of a 
system of the form \eqref{eq:lve} cannot be characterized by the time-dependent eigenvalues of $A(t)$ (see example at bottom of page 3 of \cite{Coppel}).  

The time-dependent stability spectra we consider in this paper are the Lyapunov spectrum and the Sacker-Sell spectrum. The Lyapunov Spectrum expresses the asymptotic exponential growth/decay rates of fundamental matrix solutions of \eqref{eq:lve}, and the Sacker-Sell spectrum is  defined in terms of exponential dichotomies of fundamental matrix solutions of \eqref{eq:lve} \cite{DVV3}.  Both of these spectra characterize the stability of the zero solution in the sense that the intersection of either spectra with the negative half-line $(-\infty,0)$ gives the existence of a stable or attracting direction of the equilibrium solution. Containment of either spectrum in the left half line implies asymptotic stability with exponential decay rate.  However, the Sacker-Sell spectrum is better suited to characterizing the stability of nonlinear problems without additional hypotheses on the associated linear variational equation \cite{DVV3}. 

Exact computation of the end-points of the Lyapunov and Sacker-Sell spectra  generally requires knowing the exact solution of \eqref{eq:lve}, and so most applications of these spectra involve using numerical methods for their approximation. For the methods in this paper, we use the continuous $QR$ method to change to a coordinate system in which the linear variational equation has an upper triangular coordinate system. This corresponding system is used for the numerical approximation of spectral endpoints since the method is advantageous in terms of accuracy and analysis.  

Generically, systems with an upper triangular coefficient matrix have  spectral endpoints that can be expressed as functions of the diagonal entries of the coefficient matrix. 
Define the Steklov averages of $\dot{y} = B(t)y$ for $t > 0$ and window length $H \geq 0$ as the $n$ quantities
\begin{equation}\label{eq:steklovavg}
\mu_i(t,H) = \frac{1}{H}\int_t^{t+H} B_{i,i}(\tau)d\tau, \quad i=1,\hdots,n.
\end{equation}
As a convention we take $\mu_i(t,H) = B_{i,i}(t)$  if $H= 0$ for $i=1,\hdots,n$. Let $\Sigma_{LE} = \cup_{i=1}^{n}[a_i,b_i]$ and $\Sigma_{SS} = \cup_{i=1}^{n}[\alpha_i,\beta_i]$ respectively denote the Lyapunov and Sacker-Sell spectra of $\dot{y} = B(t)y$.  The right-endpoints $\{b_i\}_{i=1}^{n}$ are referred to as \textit{upper Lyapunov exponents} or simply \textit{Lyapunov exponents}, and the left end-points $\{a_i\}_{i=1}^{n}$ are referred to as \textit{lower Lyapunov exponents}.  Under the  
assumption that $\dot{y}=B(t)y$ has an integral separation\footnote{Integral separation is a generic property of linear systems \cite{palmer} that can be changed to an upper triangular system that is integrally separated by a Lyapunov transformation. The $Q$ transformation from the $QR$ decomposition is one such Lyapunov transformation. Thus the assumption of integral separation is not an issue. It has many important consequences for stability spectra, as are discussed in depth in \cite{DVV3}. One such consequence is that integrally separated systems have distinct Lyapunov exponents.}
we have \cite{DVV3}:
$$a_i = \liminf_{t\rightarrow \infty}\mu_i(0,t), \quad b_i =  \limsup_{t\rightarrow \infty}\mu_i(0,t), \quad i=1,\hdots,d.$$
As long as $B(t)$ is bounded and continuous, we have:
$$\alpha_i = \liminf_{H\rightarrow \infty}\text{inf}_{t >0}\mu_i(t,H), \quad \beta_i = \limsup_{H\rightarrow \infty}\text{sup}_{t > 0}\mu_i(t,H), \quad i=1,\hdots,d.$$
As an alternative method for computing $a_i$ and $\alpha_i$ without using a liminf one may compute the $b_i$ and $\beta_i$ of the associated adjoint system $\dot{y}=-B(t)^{T}y$.

The continuous $Q$R method provides a way for transforming a general linear system of the form \eqref{eq:lve} to a system with an upper triangular coefficient matrix.  Let $X(t)$ be a fundamental matrix solution of \eqref{eq:lve} and let $X(t) = Q(t)R(t)$ be a QR factorization where $Q(t) \in \mathbb{R}^{n\times n}$ is orthogonal and $R(t) \in \mathbb{R}^{n\times n}$ is upper triangular with positive diagonal entries.  The matrix function $R(t)$ is a fundamental matrix solution of $\dot{y} = B(t)y$ where $B(t)= Q(t)^T A(t)Q(t) - Q(t)^T \dot{Q}(t)$ is upper triangular.  In practice, it is best to avoid forming the fundamental matrix solutions $X(t)$, which involves approximating either a potentially stiff initial value problem or a solution that is unbounded and growing exponentially. Instead we work with the coefficient matrices $A(t)$ and $B(t)$ and the orthogonal factor $Q(t)$.  This can be done since $Q(t)$ satisfies the following differential equation \cite{DVV1995}:
\begin{equation}\dot{Q}(t) = Q(t)S(Q(t),A(t)), \quad S(Q,A)_{ij} = \left\{
\begin{array}{rl}
(Q^T AQ)_{ij},& i>j\cr
0, & i=j\cr
-(Q^T A Q)_{ji},& i<j\cr
\end{array}
\right.
\label{eq:Qeqn}
\end{equation}
Thus we can form $B(t)$ by using $A(t)$ and solving a $n \times n$ differential equation for $Q(t)$ using some specified initial orthogonal factor $Q(0)=Q_0$.  

Although solving for the orthogonal factor $Q(t)$ involves solving $n^2$ differential equations (or $np$ if only the first $p$ columns of $Q(t)$ are needed), it is in many ways preferable to finding and then factoring the fundamental matrix solution $X(t)=Q(t)R(t)$. The global error for the orthogonal integration of $Q$ is bounded by a multiple of the local error and is not restricted by numerical stability considerations \cite{DVV0a}. Additionally, $Q(t)$ remains bounded while the fundamental matrix solution may become unbounded. This implies that the error in computing stability spectra and related quantities can be controlled in terms of the local accuracy and is robust even in the presence of unbounded solutions \cite{DVV4,DVV2}.  The columns of $Q(t)$ contain useful information about the directions of exponential growth and decay since for each $p \leq n$, the first $p$ columns of the matrix $Q(t)$ form an orthogonal basis for the $p$ Lyapunov vectors corresponding to the largest $p$ Lyapunov exponents \cite{DVVE1}. 

The Steklov averages (as given in equation \eqref{eq:steklovavg}) characterize the exponential growth or decay rates of a fundamental matrix solution over a given window of time.  To see this, let $X(t)$ be a fundamental matrix solution of \eqref{eq:lve} and take a $QR$ factorization $X(t) = Q(t)R(t)$. As noted above, $R(t)$ is a fundamental matrix solution of $\dot{y} = B(t)y$ where $B(t)$ is upper triangular, bounded, and continuous. Similarly consider the matrix differential equation $\dot{\Phi}=B(\tau)\Phi$, $\tau > t$, with initial value $\Phi(t;t) = I$. It has the unique upper triangular solution $\hat R(\tau;t)$.

Then, given $t > 0$ and $H > 0$, we can express $R(t+H)$ as a matrix product $R(t+H) = \hat R(t+H;t)R(t)$. Since $\hat R(\tau;t)$ is upper triangular, it can be decomposed into a matrix of diagonal entries and a strictly upper triangular nilpotent matrix $\mathcal{N}(\tau;t)$:  
$$\hat R(t+H;t) = \text{diag}\left(e^{H \mu_1^B(t,H)},\hdots,e^{H \mu_n^B(t,H)} \right) + \mathcal{N}(t+H; t) \equiv \hat R^{\text{diag}}(t+H;t) + \mathcal{N}(t+H;t).$$
Using any matrix norm with a submultiplicative property, such as an $\mathcal{L}_p$ matrix norm, $\|\mathcal{N}(\tau;t)\|\leq 1$ and 

\begin{equation}\label{eq:steklovgrowth}
\|R(t+H)\| = \|\hat R(t+H;t)R(t)\| \leq \|\hat R^{\text{diag}}(t+H;t)\|\cdot \| R(t)\|.
\end{equation}
The average exponential growth/decay rates of $\hat R^{\text{diag}}(t+H;t)$ are given by the Steklov averages of $B(t)$ of window length $H$, $\mu_i(t,H)$, $i=1,\hdots,n$. Thus equation \eqref{eq:steklovgrowth} can be interpreted as saying that the average exponential growth/decay of $R(t)$ on the interval $[t,t+H]$ is determined by the Steklov averages of $B$. Since $X = QR$ and $Q$ is orthogonal, the Steklov averages are also a measure of the average exponential growth/decay rates of $X$ on the interval $[t,t+H]$. This argument implies that Steklov averages of window length $H>0$ are a measure for the transient exponential stability (attractivity or repulsivity) of a trajectory $x(t,t_0,x_0)$ on the interval $[t,t+H]$.  Hence it is changes in Steklov averages and $Q(t)$ on intervals of the form $[t,t+H]$ as opposed to approximate Lyapunov and Sacker-Sell spectral end-points and the limiting behavior of $Q(t)$ that we use as indicators of rate-tipping.  

Examining transitions in dynamics as a critical rate changes is reminiscent of a bifurcation, however the theory of bifurcations in nonautonomous systems is still underdeveloped. Most papers on bifurcations in nonautonomous systems build towards extending definitions of Lyapunov stability using pullback attractors \cite{KS2005, LRS2002, LRS2006}. In the rate-tipping examples in this paper, the linear ramping case is possibly a nonautonomous bifurcation of pullback attractors. In \cite{rasmussen2010}, Rasmussen considers the notion of nonautonomous bifurcations on only a finite time interval and uses exponential dichotomy to define a finite-time stability spectrum with which to classify bifurcations of attracting and repelling invariant sets. We build off these ideas by using $QR$-based numerical methods for approximating time-dependent stability spectra and Steklov averages to detect rate-induced tipping.

\section{Methods}
\label{Section-Methods}

In this paper, we examine several examples in order to explore the definition of ``tipping," and specifically looking for indicators of what time a trajectory tips. We consider the definition of tipping using a tracking radius as given in \cite{ashwin2012} and two other numerical methods to determine the tipping time: the time series of Steklov averages, and the $Q$ matrix from the $QR$-method. With each of these numerical methods, we define a threshold to determine tipping within the system. Our goal is to answer the question of `will this tip soon?' not just `will this tip?

\textbf{\textit{Tracking Radius Method.}} In \cite{ashwin2012}, Peter Ashwin et al provide a sufficient criterion for rate-tipping, where a solution has not tipped if it stays within a certain tracking radius of a stable QSE. Thus a trajectory tips if it leaves that radius. However, this method determines only if a trajectory will inevitably tip, and not the time where it will tip. We define the tracking radius to be the distance between the system's stable and unstable QSE, and we calculate the time a trajectory tips by the time value at which it leaves this tracking radius. 

Since the two-dimensional Resource-Consumer model (see Table \ref{table:equations} for equations) has a stable and  saddle QSE but no unstable QSE, we do not apply the tracking radius method to this problem. If a tracking radius were defined using another distance for the Resource-Consumer model, it would detect tipping for any choice of radius $\rho<6$ (where 6 is the distance between the stable and saddle QSE's) due to the fact that the resource is decreasing away from the stable QSE to $R=0$ for any value of the rate parameter.  Above the critical rate of $r=-.002$, the tracking radius method would detect ``tipping'' in the Resource-Consumer model as simply the slow drift away from the stable QSE. At and below the critical rate of $r=-.002$, the resource exhibits a qualitatively different behavior wherein after a slow drift away from the stable QSE, the resource rapidly decreases to zero over a short window of time (see Figure \ref{fig:Resource-Consumer}).  Because the tracking radius method detects tipping whether or not the resource trajectories exhibit rapid shifts in value, we eliminate it as a possible method of analysis of this system.

\textit{\textbf{Steklov Averages Method.}} Plots of the time series of Steklov averages with window length $H=2$ accompany all of the examples. We found that changing the window length changed the scale of the Steklov averages, but not the qualitative shape of the time series. 

By observation, Steklov averages for untipped trajectories converge to a constant, negative value as $t$ gets large, or rather that for any $\epsilon>0$  there exists a time $T$ such that for any $t>T$ the absolute value of the time-derivative of the time series will be less than or equal to zero. However, we observed that the Steklov averages for a trajectory that tips and diverges  does not converge, and instead continues to increase to positive values. We define a threshold for tipping and the time at which the trajectory tips by looking at the time-derivative of the Steklov average time series, and then finding when the time-derivative first becomes positive. To avoid false indications of tipping as trajectories settle towards an attractor, we compare the Steklov average time series for two trajectories with the same initial conditions, but where one is a solution to the system where the rate is less than the critical rate and we know the trajectory does not tip, and where the other is a solution to a system with a different rate parameter and we do not know if the trajectory tips.

For $i=1,\hdots,n$ let $\mu^{B_*}_i(t,H)$ denote the Steklov averages of an untipped trajectory, and $\mu_{i}^B(t,H)$ be the Steklov averages of the trajectory in question. After identifying the time $T>0$ so that for all $t> T$
\[ \left\lvert \frac{d \mu^{B_*} _i}{dt} \right\rvert \leq \epsilon\qquad i=1,\ldots,n\]
for some small $\epsilon>0$, we consider the time rate of change of the $\mu_{i}^B(t,H)$ Steklov average time series for $t>T$. If $\mu_{i}^B$  on the interval $[t, t+h]$ satisfies
\[\left\lvert \frac{d \mu^{B_*} _i}{dt} \right\rvert> \epsilon \qquad i=1,\ldots,n\]
for some $t>T$ and some fixed $h>0$, we say that we detect tipping in the system and report the tipping time as $t+h$.  In our examples we take $\epsilon=.001$ and $h=1$.  

\textit{\textbf{Q-angle method.}} Our observation that the direction of growth for the tipped and untipped trajectories in some of our example systems were in opposite directions led us to develop the $Q$-angle method to be able to measure this change. The $Q$-angle method uses the leading Lyapunov vector, the direction at which perturbed trajectories grow/decay at the rate of the leading Lyapunov exponent, to find a change in the direction of the leading growth term so as to to detect tipping. This works only for systems with $n\geq 2$ as it compares the angle between the first column of the $Q$ matrices for two different solutions. However it is possible to apply the $Q$-angle method to a one-dimensional system by extending the system to a partially decoupled two-dimensional nonautonomous system where the first state variable $x$ has the same equation as the corresponding one-dimensional example, and the second state variable $y$ depends on $x$.

Let $Q_*(t)$ denote the time series of the $Q$ matrix from the $QR$ decomposition of a trajectory that does not tip and let $Q(t)$ denote the time series of the $Q$ matrix from the $QR$ decomposition of a trajectory with the same initial condition and different parameter rate, where we do not know if the trajectory tips. For each time $t_n$ in the time series, the $Q$-angle time series is given by $\arccos(Q_{*,1} (t_n)\cdot Q_1(t_n))$ where $Q_{*,i}(t)$ and $Q_i(t)$ denote the $i^{th}$ column of $Q_*$ and $Q$ respectively.  In the two-dimensional examples in this paper, we observed a characteristic dip then followed by an increase which seemed to be correlated with tipping. We define a threshold for tipping by the minimum of this dip, the time step of this minimum is the value we report for the Q-angle tipping time. The theoretical cause of this signature is unclear and could be the subject of further investigation. 

\textbf{\textit{Implementation of Numerical Methods.}} We approximate the Steklov averages and $Q(t)$ by first fixing an initial condition $x(0) = x_0$ and initial orthogonal $Q(t_0) = Q_0$.  We approximate $Q(t)$ by using a standard Matlab ODE solver to integrate equations \eqref{eq:ode} and \eqref{eq:Qeqn} simultaneously, while being careful to project $Q(t)$ at each step to maintain its orthogonality, and then forming $A(t)$ as the linearization around the approximate solution trajectory.  We then use the approximate $Q(t)$ and $A(t)$ to form $B(t)$ from which we approximate the Steklov averages $\mu_i^B(t,H)$ with a numerical quadrature rule.

\FloatBarrier
\section{Examples}
\label{Section-Examples}

\begin{table}[h] 
\begin{center}
\addtolength{\tabcolsep}{1mm}
\renewcommand{\arraystretch}{1.2}
\resizebox{\textwidth}{!}{%
\begin{tabular}{||p{1.5cm}|p{6cm}|p{2cm}|p{1.4cm}|c|c||}
\hline
\textbf{Problem}  & \textbf{Equations} & \textbf{Initial Cond.} & \textbf {Time Interval} & \textbf{No Tip} & \textbf{Tip}  \\
\hline
\hline
Unique Linear     & $\dot x = - (x-\lambda(rt))(x - \lambda(rt) - \delta) $  & $x(0)=.5$ & $[0,60]$& $r< .0625$ & $r> .0625$  \\
\hline
Bistable Linear    & $\dot{x}  = - (x-\lambda(rt))((x - \lambda(rt))^2 - \delta^2)$ & $x(0)=.5$ & $[0,140]$ &$r<.048$ & $r\geq .049$ \\
\hline 
Unique Logistic  & $\dot x = - (x-\lambda(rt))(x - \lambda(rt) - \delta) $ & $x(0)=.5$ & $[0,60]$ & $r< .5$     & $r> .5$          \\
\hline 
Bistable Logistic & $\dot{x}  = - (x-\lambda(rt))((x - \lambda(rt))^2 - \delta^2)$ & $x(0)=.5$ & $[0,100]$ &$r< .377$ & $r\geq .378$  \\
\hline
\hline
Bistable Linear, 2D          & $\dot{x}  = - (x-\lambda(rt))((x - \lambda(rt))^2 - \delta^2)$ $\dot y = x-y$ & $x(0)=.5$  $y(0)=.5$ & $[0,140]$  & $r<.048$ & $r\geq .049$ \\
\hline 
Bistable Logistic, 2D        & $\dot{x}  = - (x-\lambda(rt))((x - \lambda(rt))^2 - \delta^2)$ $\dot y = \frac{1}{2}x^2-y$  & $x(0)=.5$  $y(0)=.5$ & $[0,100]$ & $r< .377$ & $r\geq .378$ \\
\hline
Resource-Consumer Linear    & $\dot R = \lambda(rt)R\left(1-\frac{R}{K}\right) -\frac{aCR}{R+R_h} $ \hspace{10mm} $\dot C = \epsilon\left(\frac{eaCR}{R+R_h}-mC\right)$ & $R(0)=6$ $C(0)=16$ & $[0,2500]$ & $ r>-.001$ & $r\leq -.002$\\
\hline
\end{tabular}}
\caption{Equations, initial conditions, time interval of solution, and tipping values for the example systems considered in this paper.  The one dimensional problems and the two dimensional bistable problems were solved using MATLAB's \texttt{ode45} with \texttt{`RelTol'} and \texttt{`AbsTol'} set to \texttt{1e-14}. In all problems $\delta=.5$.  The Resource-Consumer model \cite{siteur2016ecosystems} was solved using MATLAB's \texttt{ode15s} with \texttt{'RelTol'}$=0.6$, \texttt{'AbsTol'}$=0.6$,\texttt{'BDF'} set to \texttt{'on'}, \texttt{'MaxOrder'}$=1$, \texttt{'MaxStep'=1e-3}; and parameters $a=1$, $e=1$, $K=10$, $m=.75$, $\epsilon=.01$, and $R_h=2$.} \label{table:equations}
\end{center}
\end{table}

Table \ref{table:equations} provides the details of the nonautonomous examples considered in this paper: four one-dimensional examples and three two-dimensional examples. Six of these examples were chosen as test cases for the heuristic numerical methods we propose based on their relationship to examples in \cite{ashwin2012, ashwin2015,perryman2014}, and the seventh is the Resource-Consumer model in \cite{siteur2016ecosystems}. All examples depend on a parameter $\lambda(rt)$ that is an explicit function of time. Examples labeled with `Unique' in the title have one stable QSE, and those labeled with `Bistable' in the title have two stable QSEs. We examine two types of nonautonomous parameters, linear ramping and asymptotically constant, in order to make direct comparison with rate-tipping examples \cite{ashwin2012, ashwin2015,perryman2014}. It is worth noting that in looking at asymptotically constant parameter drift, in contrast to the hyperbolic tangent function used by Ritchie and Sieber (also used by \cite{ashwin2012,perryman2014,ashwin2015}), we use a parameter that satisfies a logistic differential equation. The labels `Logistic' and `Linear' distinguish the two cases of nonautonomous parameter $\lambda(rt)$.  Problems with `Logistic' in the title have a parameter $\lambda(t)=(e^{r t})/(\lambda_1+e^{r t})$ where $\lambda_1=(1-\lambda_0)/(\lambda_0)$ and $\lambda_0=10^{-6}$, so that $\lambda(t)$ satisfies the logistic differential equation
\begin{equation}
\label{EQ-lambda-dot-log}
\dot{\lambda}(t)  = r \lambda(t) (1-\lambda(t)).
\end{equation}
Problems with `Linear' in the title mean that $\lambda(t)=rt+\lambda_0$
where we take $\lambda_0=0$. Here $\lambda(t)$ satisfies the constant differential equation
\begin{equation}
\label{EQ-lambda-dot-lin}
\dot{\lambda}(t)  = r .
\end{equation}

For each example we calculate the critical tipping rate for the trajectory starting with $x(0)=0.5$. These are contained in Table \ref{table:equations}. Table \ref{table:tipping} give the results for the various methods we investigated. These results are discussed in further detail in the next two sections.

\begin{table} 
\begin{center}
\addtolength{\tabcolsep}{1mm}
\renewcommand{\arraystretch}{1.2}
\begin{tabular}{||p{1.5cm}||p{1.7cm}|p{2cm}||p{2cm}|p{2cm}||p{2cm}||}
\hline
\textbf{Problem}  & \textbf{Tracking Radius} & \textbf{Radius Tipping Time} & \textbf{Steklov Average Window Length} & \textbf{Steklov Average Tipping Time} & \textbf{Q-angle Tipping Time} \\
\hline
\hline
Unique Linear     & $.5$ & $53.94$ & $H=2$ & $38.14$ & N/A\\
\hline
Bistable Linear   & $.5$ & $ 104.9$ & $H=2$ & $36.06$ & N/A\\
\hline 
Unique Logistic    & $.5$ & $37.91$ & $H=2$ & $36.19$  &  N/A\\
\hline 
Bistable Logistic & $.5$ & $47.77$ & $H=2$ & inconclusive  & N/A\\
\hline
\hline
Bistable Linear, 2D  & $.5$ & $ 104.9$ & $H=2$ & $92.87$  & $109.8$\\
\hline 
Bistable Logistic, 2D  & $.5$ & $47.77$ & $H=2$ &  inconclusive  & $69$\\
\hline
Resource-Consumer Linear  &  N/A & N/A & $H=2$ & inconclusive & $1589$\\
\hline
\end{tabular}
\caption{Time to detect tipping using the Tracking Radius method, Steklov Average method, and Q-angle method for the minimum value of $r$ that causes tipping in the system within the time interval given in Table \ref{table:equations}. For methods where a comparison is needed, the critical values (i.e. $r=r_c$) are not used, and a value below $r_c$ is used.} \label{table:tipping}
\end{center}
\end{table}

\FloatBarrier
\subsection{One-dimensional examples}

In all four one-dimensional examples, it is possible to analytically determine the critical rate at which our given trajectory tips by considering them as two-dimensional autonomous systems where $\lambda$ is a state variable. In case where the parameter is coordinate shift $x-\lambda$, as in our examples, and in the case where that parameter is changing at a linear rate $r$, there exist invariant lines where $x$ is also changing at the same rate $\dot{x}=r$. Tipping occurs when $|r|$ is large enough that invariant lines collide and annihilate. This process could also be considered a saddle node bifurcation of pullback attactors in the one-dimensional nonautonomous system. In the Logistic examples, the critical rate at which the given trajectory tips corresponds to the extended two-dimensional autonomous system passing through a unique heteroclinic connection between two hyperbolic equilibria. Using this theory we are able to find the critical tipping rate and calculate the time at which the trajectory leaves the critical tracking radius. We compare this time to the time at which the Steklov average method indicates tipping.

\begin{figure}[h]
        \centering
         \includegraphics[width=.455\textwidth]{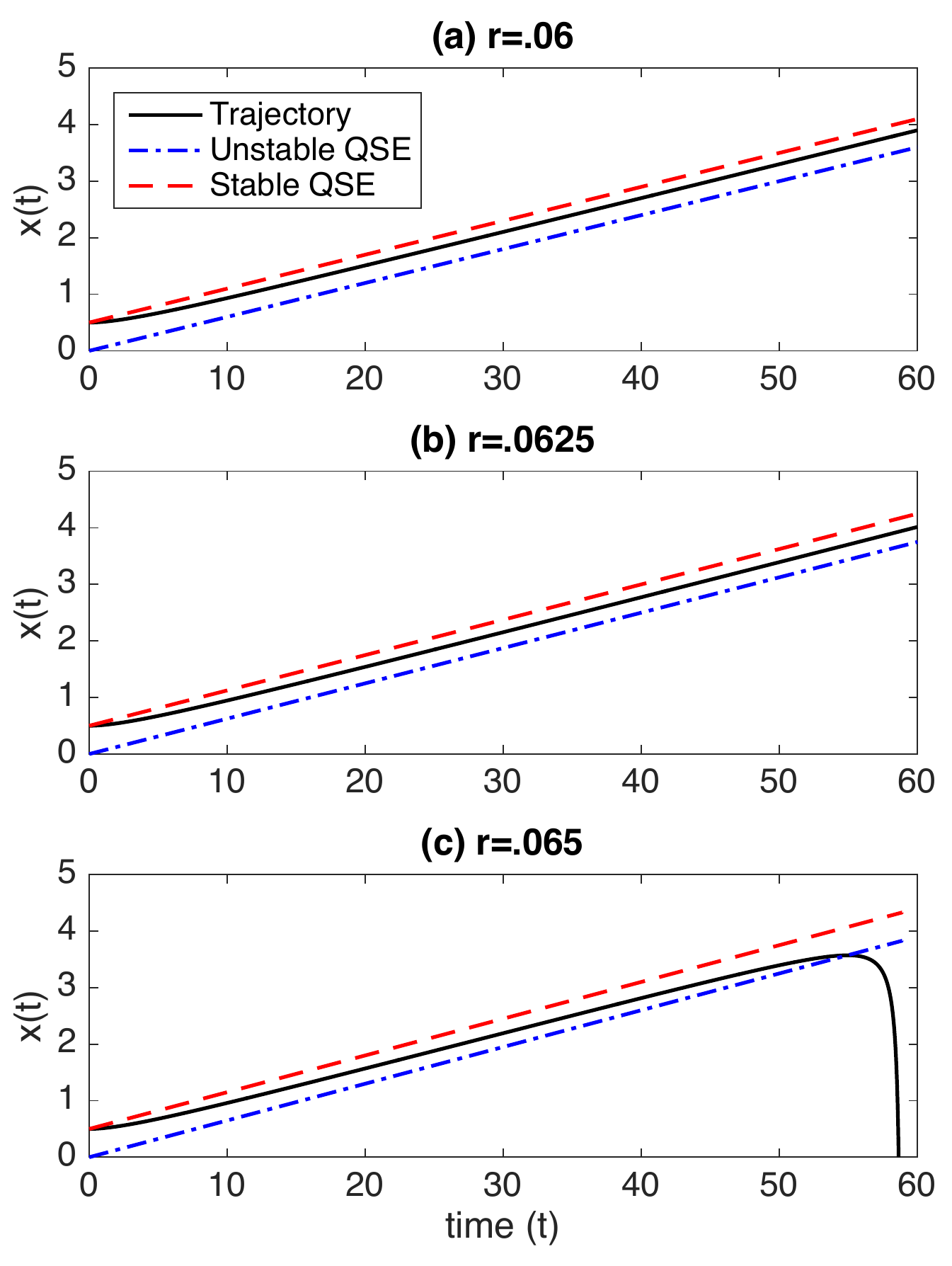}
         \includegraphics[width=.45\textwidth]{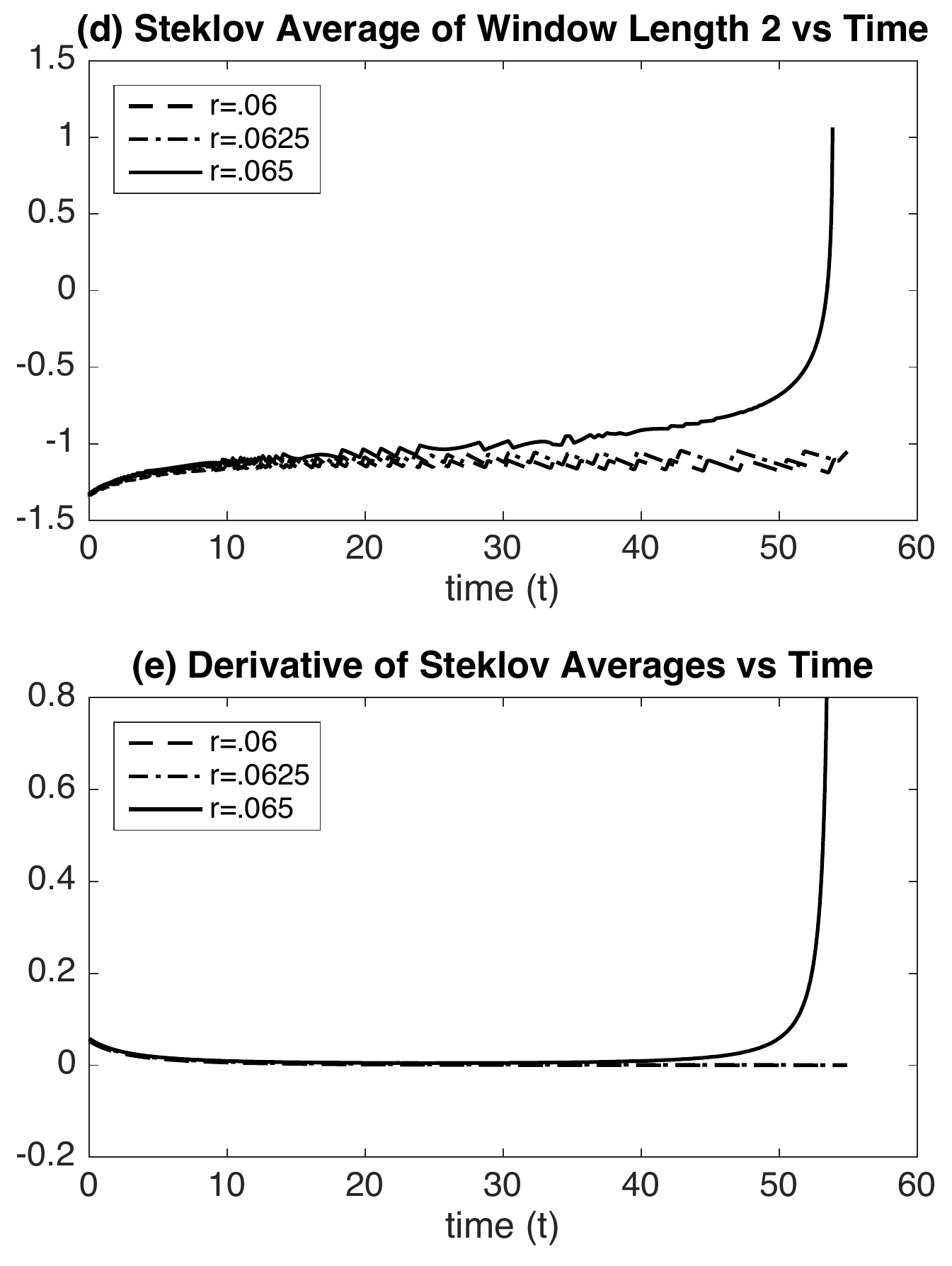}
        \caption{Rate tipping in a system with a unique stable QSE under  linear parameter growth.  In the left column are plots of trajectories (black) for different rates, $r$, along with the stable and unstable QSEs. In the right column are the time series the Steklov averages and their derivatives. See text for a detailed description.}
        \label{fig:UnLin}      
\end{figure}

\textit{\textbf{Unique Linear Example.}} For the Unique linear example, the invariant sets of the corresponding two-dimensional nonautonomous system are the solutions to 
\[\dot{x}=(x-\lambda)(x-\lambda-\delta)=r\]
These two invariant lines exists for $r<.0625$ and collide when $r=.0625$. There are no invariant lines for $r>.0625$. Therefore we consider the critical tipping rate to be $r=.0625$ and find that the given trajectory tips for $r\geq 0.0625$. In this example, the trajectory with $r=.065$ leaves the tracking radius around the stable QSE at time $t=53.94$.  By using our Steklov averages method we predict this tipping at $t=36.29$. See Figure \ref{fig:UnLin} to see the extended phase space of the system as well as time series for the Steklov averages and their derivatives.

\begin{figure}[h]
     \begin{center}
     \includegraphics[width=.46\textwidth]{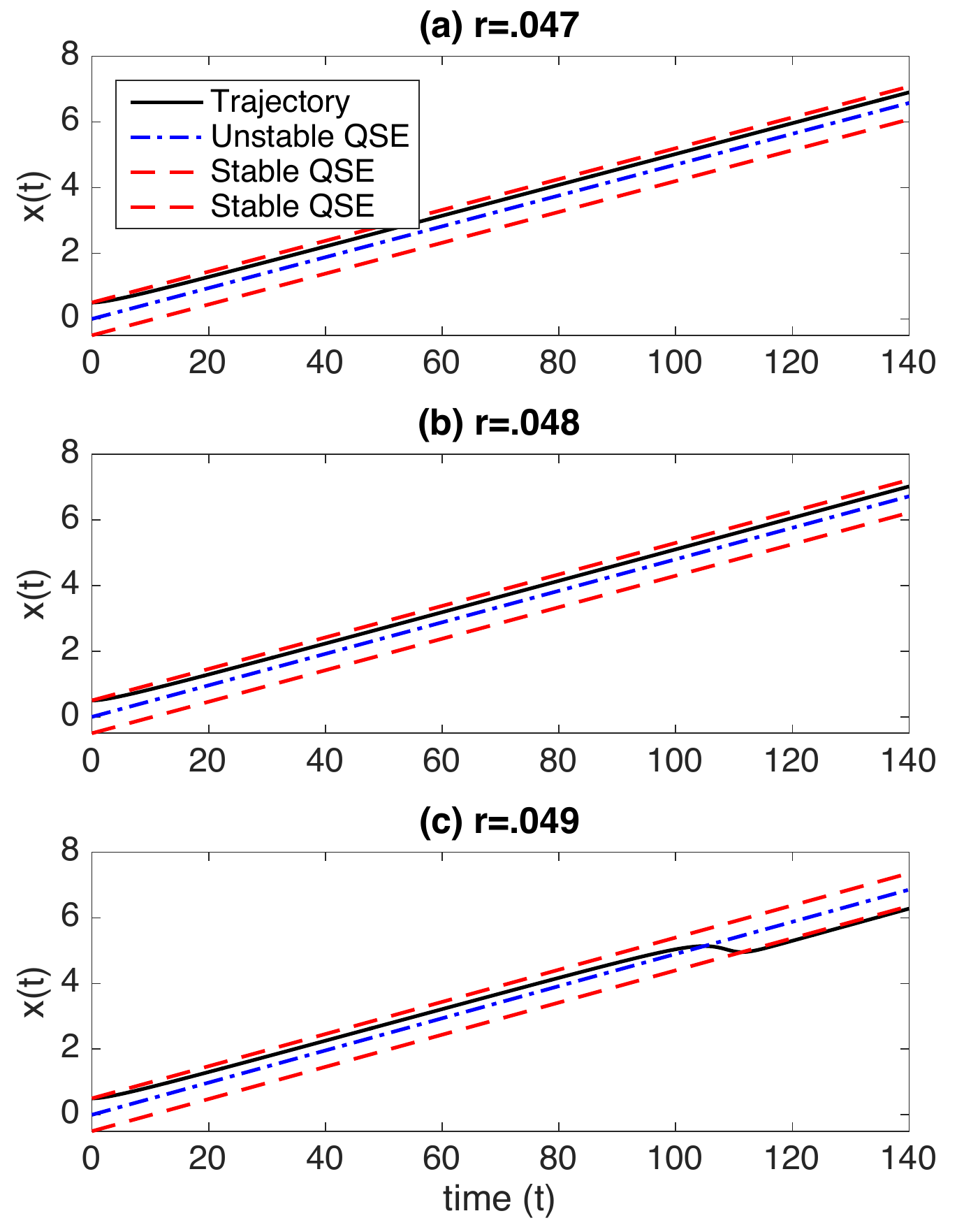}
     \includegraphics[width=.46\textwidth]{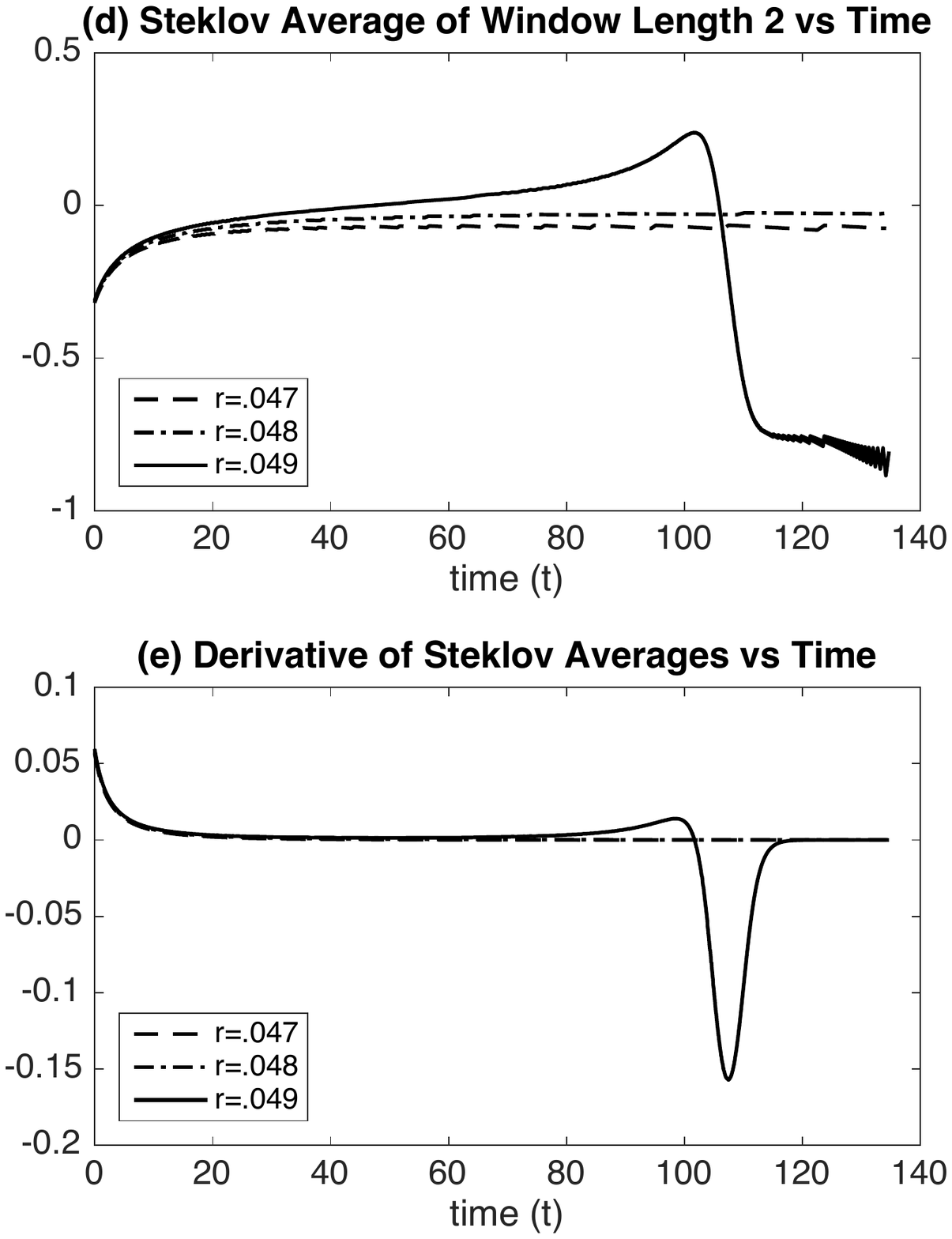}
     \end{center}
     \caption{Rate tipping in a system with a two stable QSEs under linear parameter growth. In the left column are graphs of trajectories (black) for different rates, $r$, along with two stable and one unstable QSEs. In the right column are the time series for the Steklov averages and their derivatives. See text for a detailed description.}
     \label{fig:BiLin}
\end{figure}

\textit{\textbf{Bistable Linear Example.}} For the bistable linear example, the invariant lines in the corresponding two-dimensional autonomous system are solutions to the cubic polynomial
\[\dot{x}=(x-\lambda)(x-\lambda-\delta)(x-\lambda+\delta)=r.\]
The critical tipping rate corresponds to when two invariant lines collide in phase space, and for $r\geq.049$ there is only one invariant line. In this example, at $r=.049$, the trajectory leaves the tracking radius of 0.5 at time $t=104.9$. By using our Steklov averages method we predict this tipping at $t=92.87$. Graphs of the solution of the system, Steklov averages, and derivatives of Steklov averages are in Figure \ref{fig:BiLin}.

\begin{figure}[b]
        \centering
        \includegraphics[width=.46\textwidth]{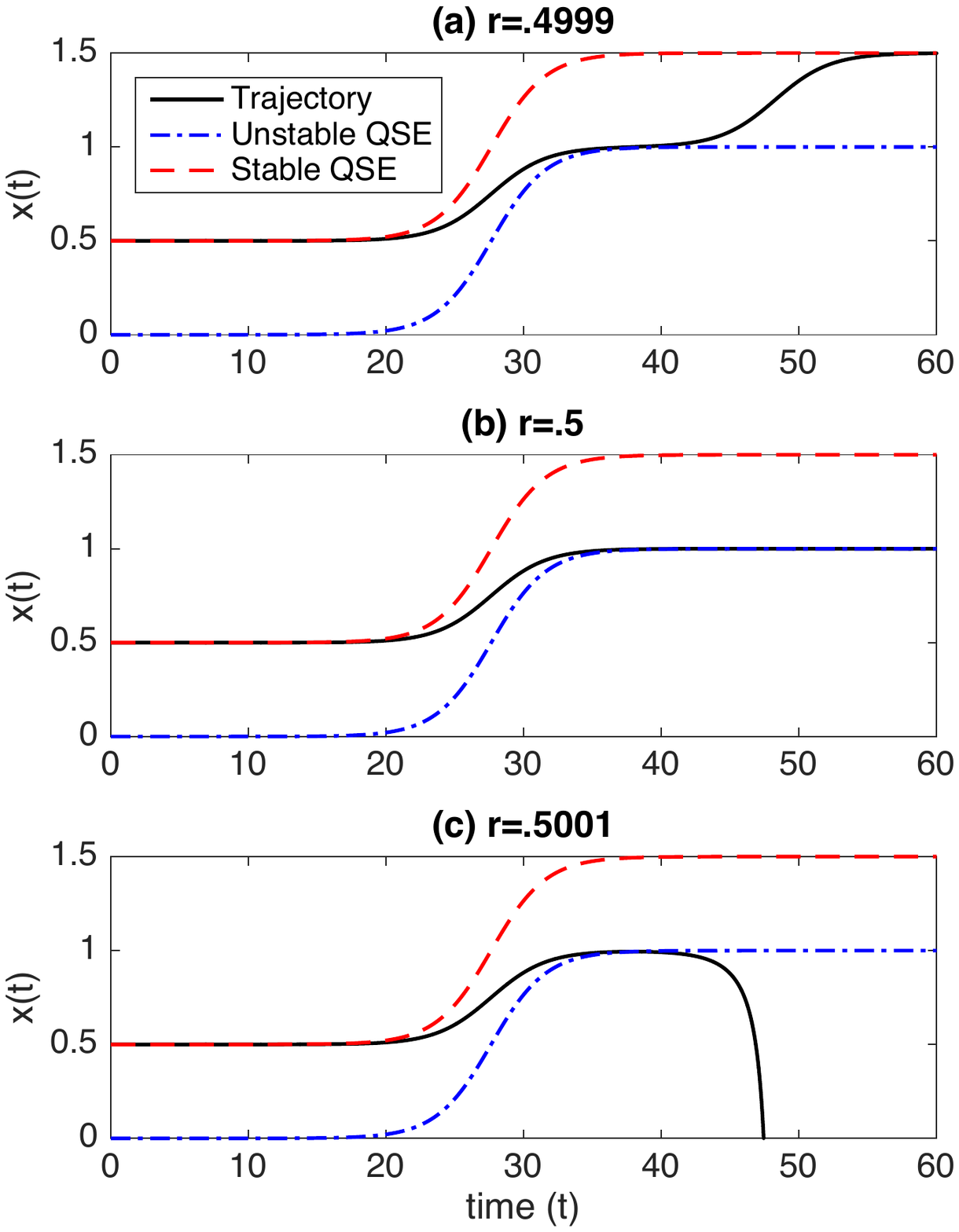}
         \includegraphics[width=.44\textwidth]{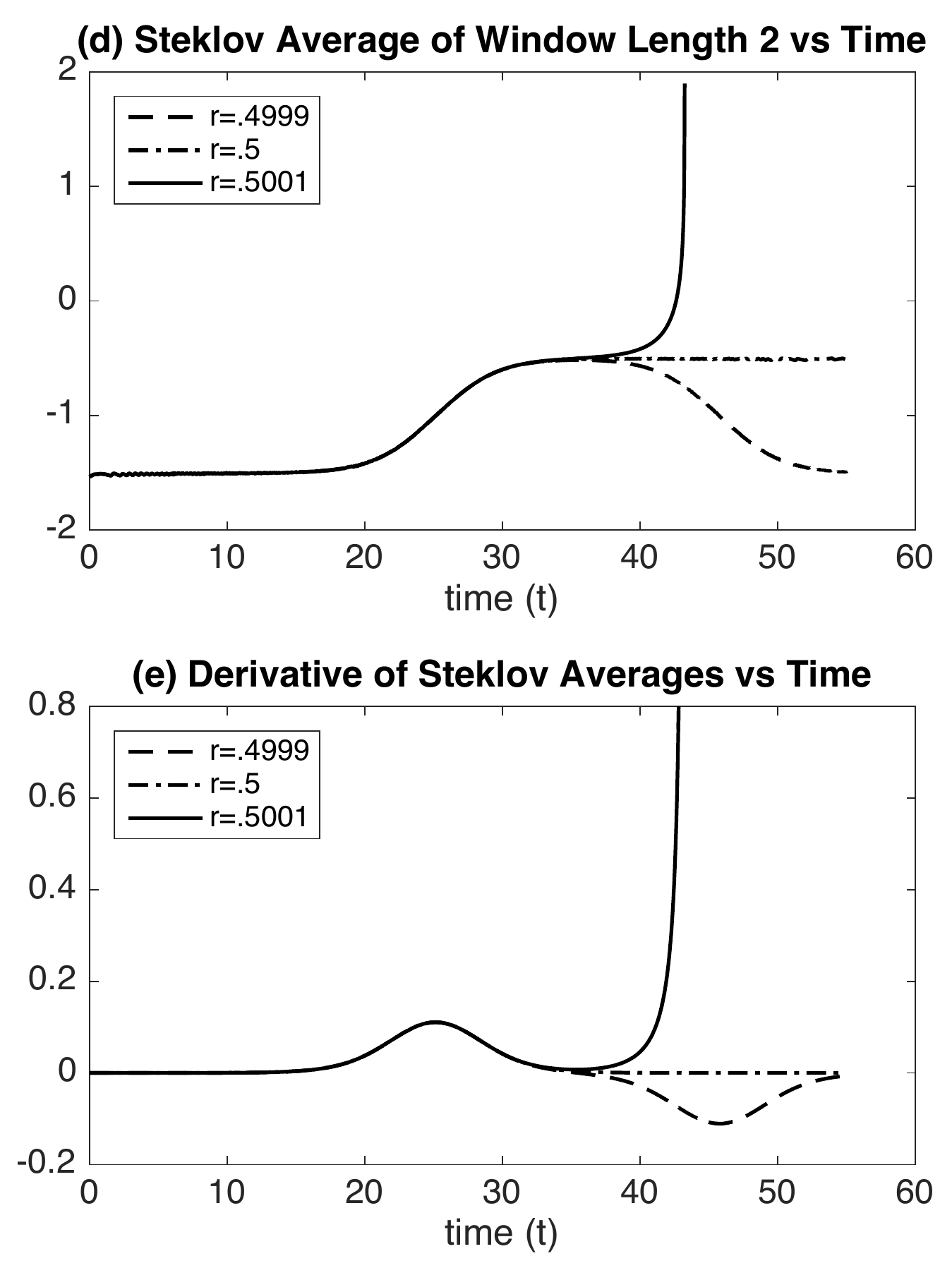}
        \caption{Rate tipping in a system with a unique stable QSE under asymptotically constant parameter growth. On the left are trajectories (black) for different rates, $r$.  On the right are the stability spectra for the trajectories under consideration.  See text for a detailed description.}
        \label{fig:UnLog}      
\end{figure}

\textit{\textbf{Unique Logistic Example.}} The Unique Logistic example has asymptotically constant parameter drift with one stable QSE. When this example is construed as a two-dimensional autonomous system, the heteroclinic connection occurs for the unique value $r=.5$. The existence of the heteroclinic connection can be shown via a Melnikov integral calculation.\footnote{The authors can provide such calculations for inclusion in an appendix, although the computations are relatively straightforward.} The value $r=.5$ is the critical rate for tipping, as the trajectory tips for $r>0.5$. In this example, the tipping trajectory for $r=.5001$ leaves the tracking radius of .5 at $t=37.91$.   Using the Steklov average method, tipping was detected in the system at $t=38.53$. Graphs of the solutions and QSEs, as well as time series for the Steklov averages and the derivatives are shown in Figure \ref{fig:UnLog}.

\begin{figure}[b]
      \begin{center}
      \includegraphics[width=.47\textwidth]{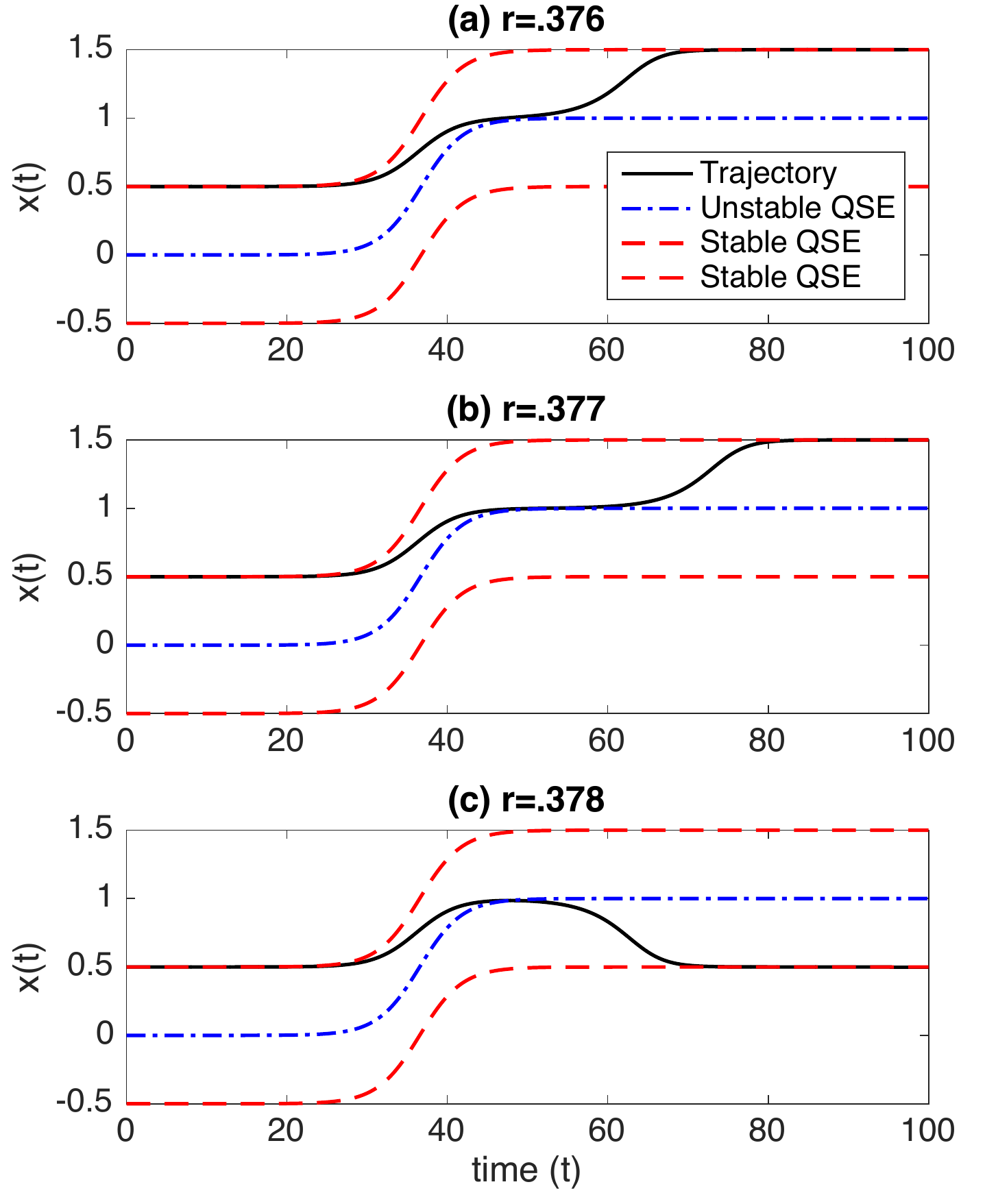}
      \includegraphics[width=.46\textwidth]{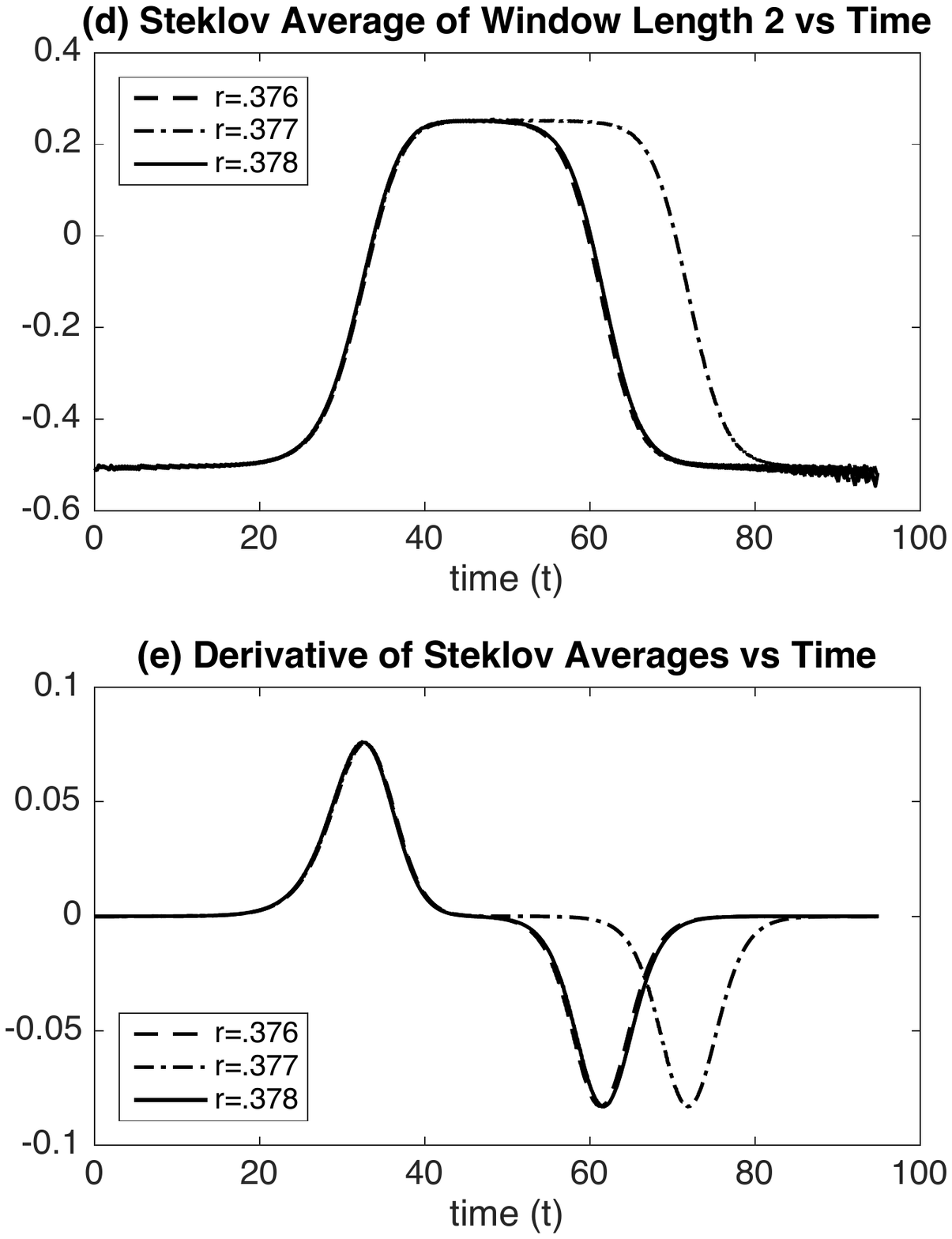}
      \end{center}
      \caption{Rate tipping in a system with two stable QSEs under asymptotically constant parameter growth. On the left are trajectories (black) for different rates, $r$.  On the right are the stability spectra for the trajectories under consideration. See text for detailed description.}
      \label{fig:BiLog}
\end{figure}

\textit{\textbf{Bistable Logistic Example.}} The last one dimensional example we consider is the Bistable Logistic, which has asymptotically constant parameter drift and two stable QSEs. This problem can be considered as a two dimensional autonomous system which has a heteroclinic connection, which numerically we found to be at $r=.377$ up to three decimal places. Thus the critical rate of tipping is approximately $r=.377$.  In this example, at $r=.378$, the trajectory leaves the tracking radius at time $t=47.77$. This example is interesting as the Steklov average method was unable to detect tipping. In the case with two stable QSE, there are solutions which exhibit \textit{endpoint tracking} (see \cite{perryman2014}), where the solution asymptotically approaches the same stable QSE as $t\to\pm\infty$. Additionally, there are solutions which tip and track the other stable QSE. In both situations, these trajectories seem to track the unstable QSE for a short amount of time. This behavior is shown in the first and third graphs in Figure \ref{fig:BiLog}. 
As the tipped trajectory then tracks the other stable QSE instead of diverging to negative infinity, the growth rates detected by the Steklov averages do not increase from negative to positive in a way that indicates instability. The inability of the Steklov averages method to detect tipping in this example motivated us to consider examples in two dimensions.

\FloatBarrier

\subsection{Two dimensional examples}

\textit{\textbf{Bistable linear and logistic examples.}} In these two examples, we have extended the one-dimenisonal examples with the same name to a partially decoupled nonautonomous system with the same $x$ differential equation and $y$ differential equation depending on $x$. Because $x(t)$ is unchanged from the one dimensional cases, the Steklov averages and the time to cross over the unstable QSE remain unchanged as well. Using the Q-angle method we were able to detect tipping in the system at $t=109.8$ and $t=69$ for the linear and logistic bistable problems, respectively.

\begin{figure}[h]
      \begin{center}
      \includegraphics[width=.32\textwidth]{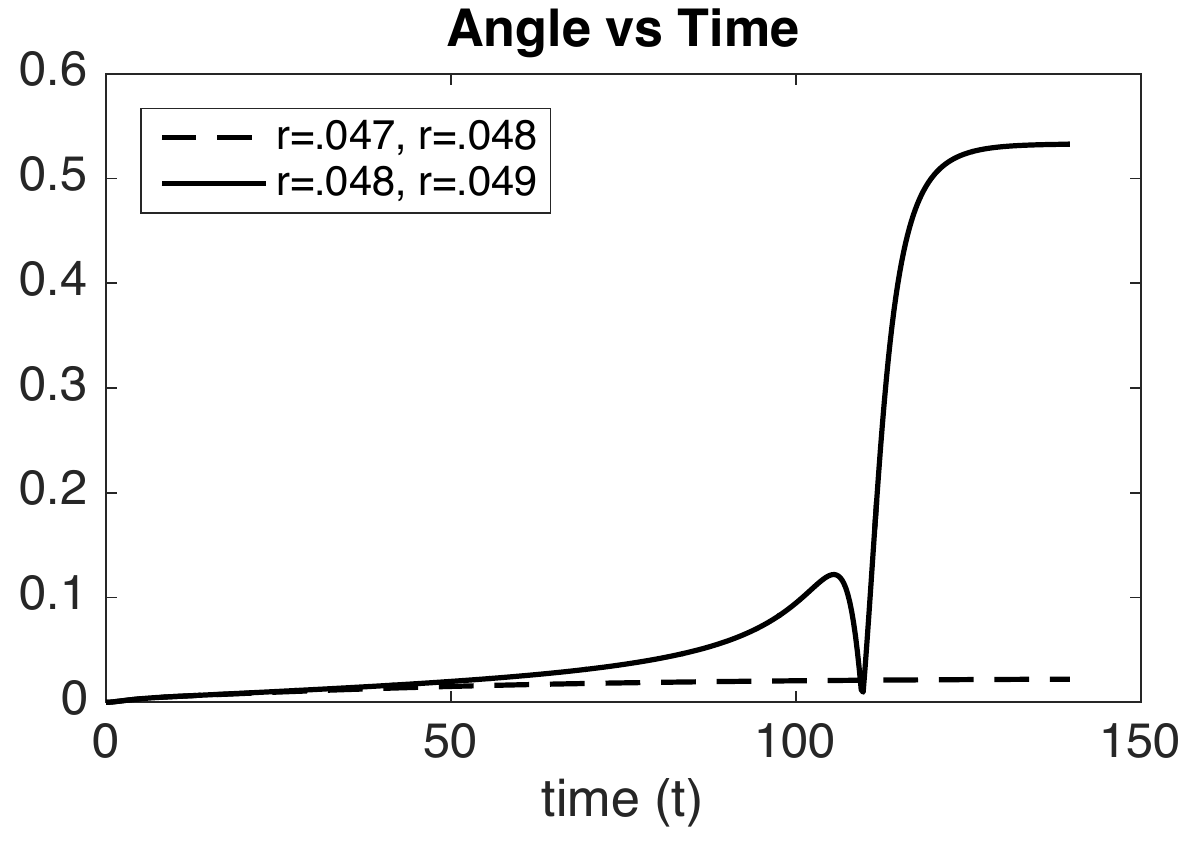}
      \includegraphics[width=.32\textwidth]{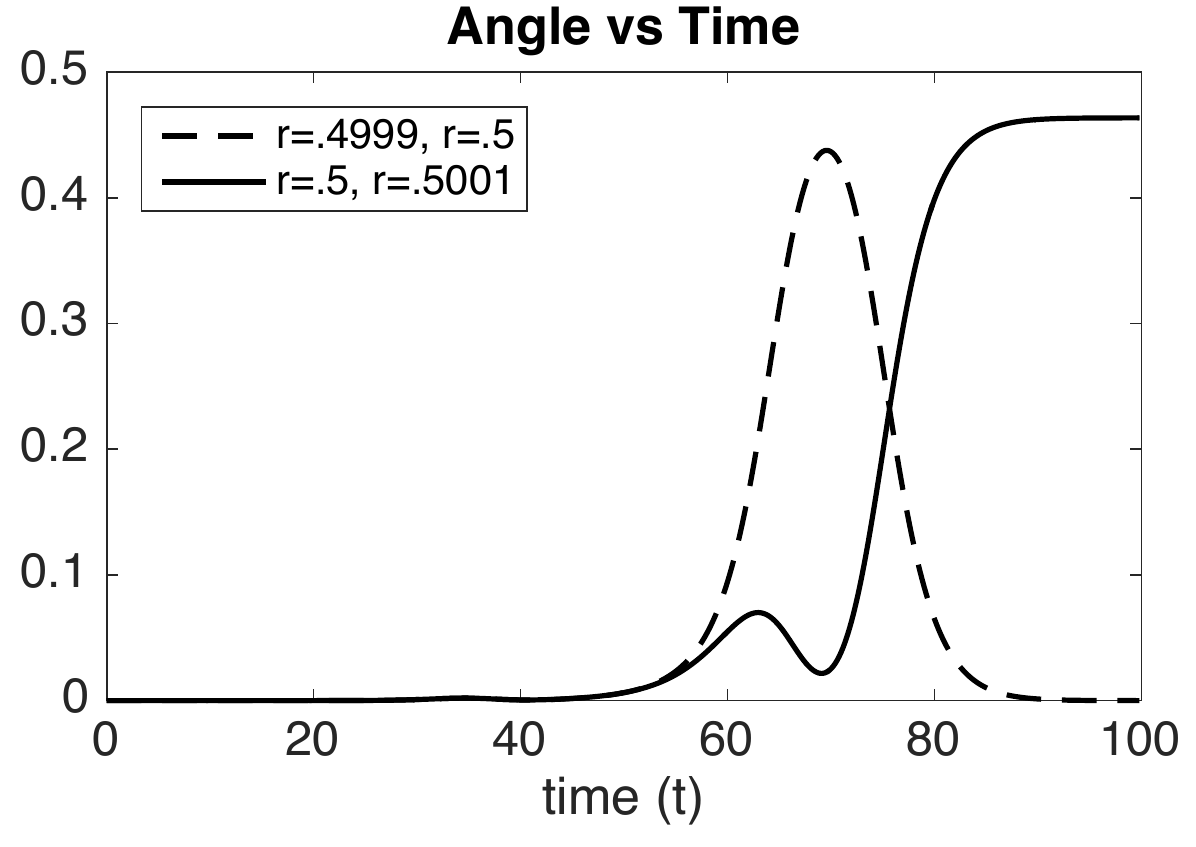}
      \includegraphics[width=.335\textwidth]{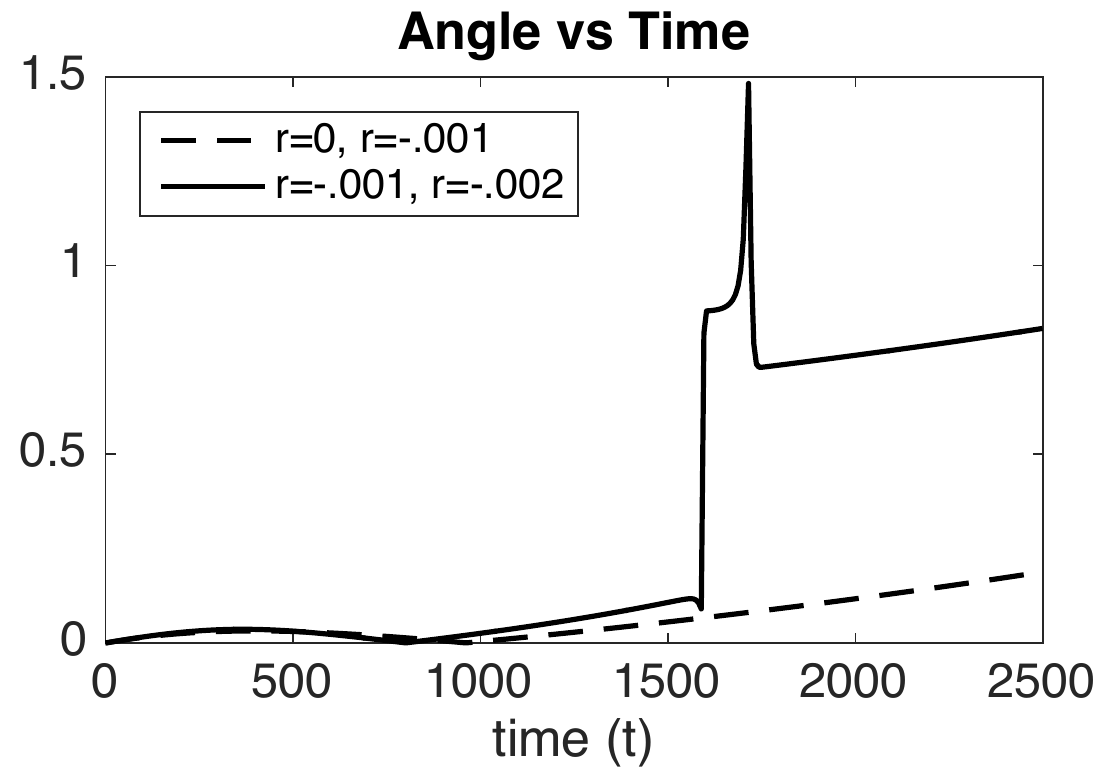}
      \end{center}
      \caption{Q angles for the Bistable Linear (left), Bistable Logistic (middle), and Resource-Consumer (right) two dimensional examples.}
      \label{fig:Resource-Consumer}
\end{figure}

\textbf{\textit{Resource-Consumer Model.}} We consider this two-dimensional example to test the robustness of the methods described in this paper. In \cite{siteur2016ecosystems} the authors give an ecological model of a resource, $R$ and a consumer of the resource, $C$, which experiences rate tipping for $r<0$.  This model has a multiplicative parameter with a linear ramping rate, as opposed to all the previous examples where the nonautonomous parameter is a linear coordinate shift. Furthermore, the ramped parameter is given by
\[\lambda(rt)=5+rt\]
where $r<0$, so contrary to our other examples, $\lambda$ is decreasing in time.  Another difference is that the ramped parameter causes the trajectory to fall away from the only stable QSE in the system and to be `attracted' to a hyperbolic QSE for a period of time. This hyperbolic QSE is also an equilibrium solution of the nonautonomous system.

The Resource-Consumer model is an interesting rate-tipping problem because it challenges the definition of tipping as not tracking a stable QSE. The given trajectory will move away from the stable QSE at $R=6$ toward $R=0$ for all values of $r<0$ and thus will inevitably tip.  Exploring the transient dynamics of the nonautonomous system shows that for $|r|$ large enough, the trajectory of $R(t)$ undergoes a quick transition towards zero.  We are able to capture this transition with the Q-angle method, but not with the Steklov averages method.

\begin{figure}
      \begin{center}
      \includegraphics[width=.53\textwidth]{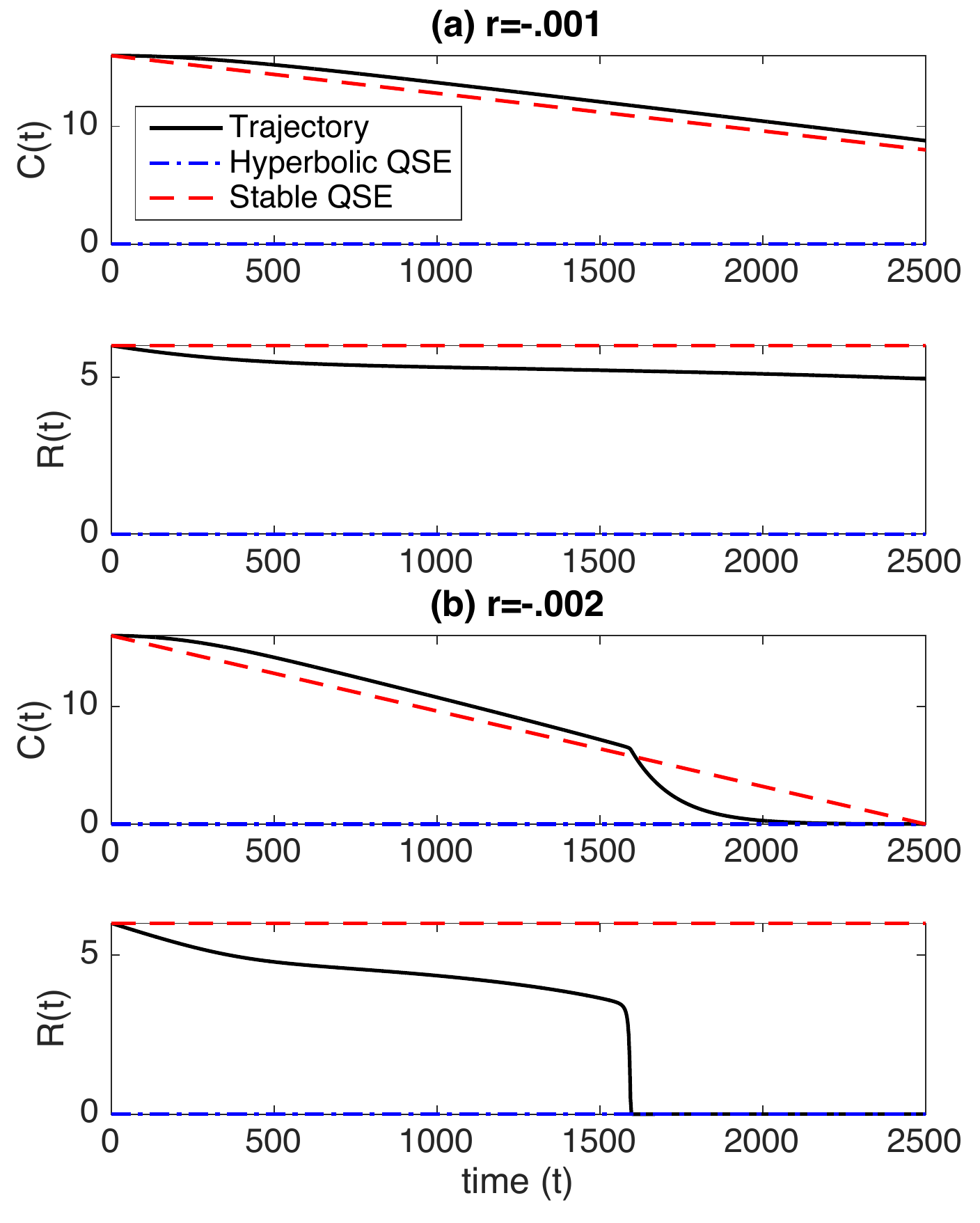}
      \includegraphics[width=.435\textwidth]{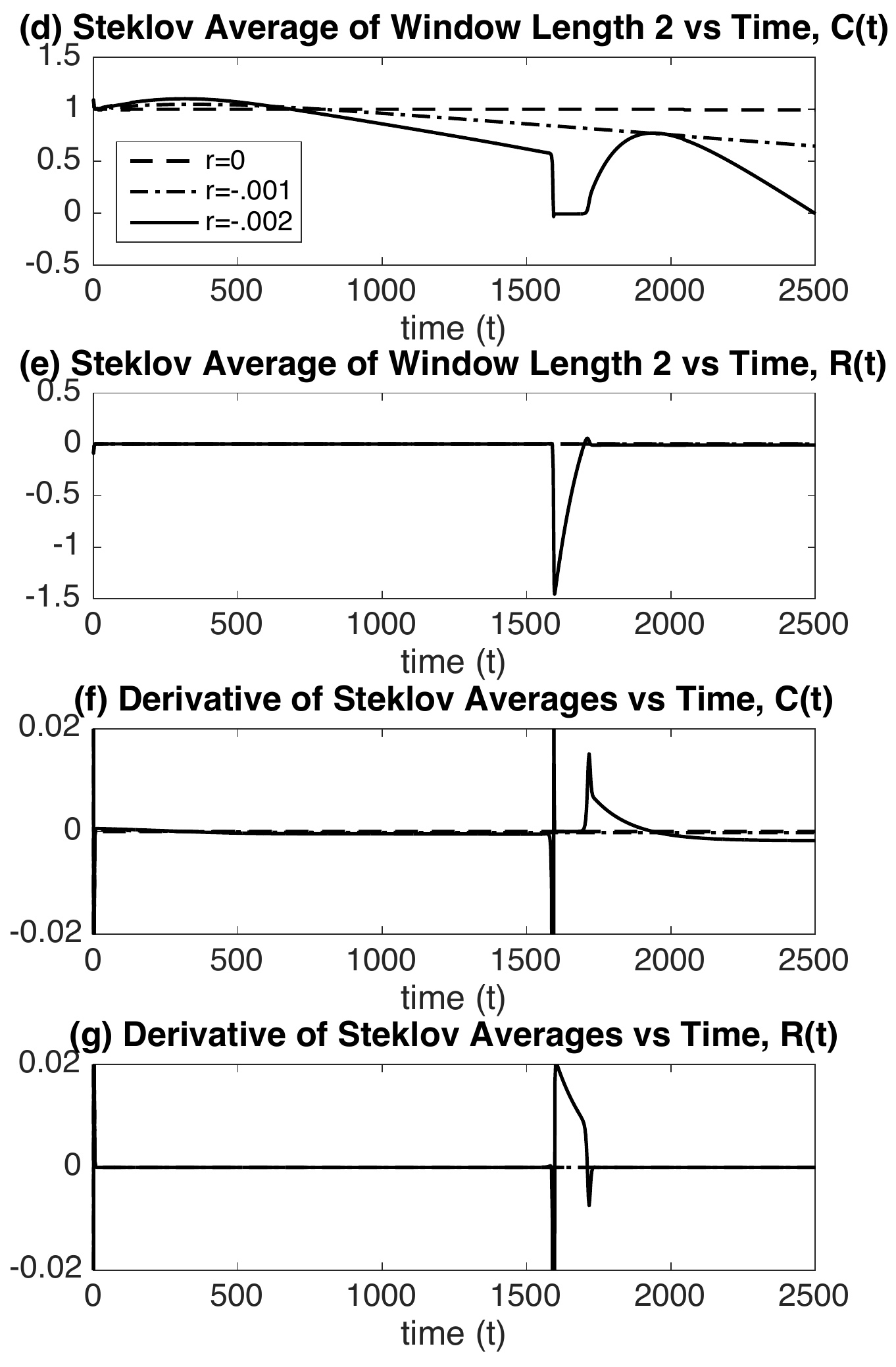}
      \end{center}
      \caption{Rate tipping in the Resource Consumer Model Rate. In the left column are trajectories (black) for different rates, $r$, for each state variable.  In the right column are the stability spectra for the trajectories under consideration. See text for detailed description.}
      \label{fig:Resource-Consumer}
\end{figure}

\FloatBarrier

\section{Discussion}
 
The mathematical study of dynamical systems has traditionally been focused on the long term behavior of solutions; however, in order to apply the study of dynamical systems to real world systems, the transient dynamics of the system must also be understood.  Understanding these transient dynamics becomes all the more important when the long term behavior is uncertain such as in the climate system.  Rate tipping and nonautonomous stability theory can be techniques which allow mathematicians to study the transient behavior of systems, but the focus must change toward transient behavior as opposed to long-term behavior of the systems in question.  In particular, the confluence of the ideas of `tracking' and `not tipping' should be avoided. The techniques outlined in this paper are intended to bring mathematicians closer to answering `Will this system tip soon?' as opposed to the question `Will this system ever tip?' Although our results are for simple systems and are not rigorously proven, we provide them here as a stepping stone toward better detection---or even prediction---of tipping points.

We found Steklov averages method indicated tipping in some cases.  In particular, the Steklov averages method worked extremely well in the one-dimensional linear ramping examples, \emph{predicting} tipping (defined as as not tracking a stable QSE) in both cases. For the asymptotically constant parameter drift example, it was only able to detect tipping for the case with one stable QSE and only indicated tipping after the trajectory had left the tracking radius. Increasing the dimension did not seem to affect the results of the Steklov averages method. For our two-dimensional examples, the Q-angle method indicated but did not predict tipping. It is unclear to us the reason why there is a characteristic dip in the Q-angle plots close in time to where the trajectory tips. This is something interesting that needs to be explored further, and we encourage the nonautonomous stability and tipping communities to look into it.

The authors acknowledge that of the many questions that were raised in this paper, very few were answered.  We hope that by providing brief backgrounds to both rate-tipping and nonautonomous stability that we may inspire others to explore the intersection of these fields further.

\section*{Acknowledgements}
The authors would like to thank Mary Silber and Eric Van Vleck for their guidance throughout this project. We also thank Sebastian Wieczorek for his suggestion to consider problems with two stable QSEs, and Kate Meyer for her suggestion to look at the Resource-Consumer model. All authors would like to acknowledge support by the AMS Mathematics Research Communities (National Science Foundation Grant No.1321794) and the Mathematics and Climate Research Network (NSF Grants DMS-0940366 and DMS-0940363). Andrew Steyer was supported by NSF Grant DMS-1419047 during this work.

\bibliographystyle{siam}
\bibliography{sources}
\end{document}